\documentclass[oneside,12pt]{amsart}

\usepackage{amscd,amsmath,latexsym,bm,pdfpages,verbatim,amsthm}
\usepackage[mathcal]{euscript}
\usepackage[normalem]{ulem}
\usepackage{graphics, setspace}
\usepackage{amssymb}          
\usepackage{epstopdf}
\usepackage[all]{xy}

\newtheorem{thm}{Theorem}[section]
\newtheorem{cor}[thm]{Corollary}
\newtheorem{lem}[thm]{Lemma}
\newtheorem{prop}[thm]{Proposition}

\theoremstyle{definition}
\newtheorem{defin}[thm]{Definition}
\theoremstyle{definition}

\theoremstyle{definition}
\newtheorem{exm}[thm]{Example}

\newtheorem{remark}[thm]{Remark}

\theoremstyle{remark}
\newtheorem*{rem}{Remark}

\newcommand{\mathsym}[1]{{}}
\newcommand{\unicode}[1]{{}}

\begin{document}
\def\X#1#2{r(v^{#2}\ds{\prod_{i \in #1}}{x_{i}})}
\def\skp#1{\vskip#1cm\relax}
\def\block{\rule{2.4mm}{2.4mm}}
\def\nd{\noindent}
\def\cxa{(\underline{CX},\underline{X})}
\def\xa{(\underline{X},\underline{A})}
\def\yb{(\underline{Y},\underline{B})}
\def\apoint{$Z(K;(\underline{X}, \underline{\ast}))$}

\def\becomes{\colon\hspace{-2,5mm}=}
\def\ds{\displaystyle}
\def\red{\color{red}}
\def\blue{\color{blue}}
\def\black{\color{black}}
\definecolor{ao}{rgb}{0.0, 0.0, 1.0}
\definecolor{beaublue}{rgb}{0.74, 0.83, 0.9}
\def\s{\sigma}
\numberwithin{equation}{section}
\title[Polyhedral Products]{Polyhedral products and features of their homotopy theory}

\skp{0.2}

\author[A.~Bahri]{A.~Bahri}\thanks{Supported in part by grant 426160 from Simons Foundation.}
\address{Department of Mathematics,
Rider University, Lawrenceville, NJ 08648, U.S.A.}
\email{bahri@rider.edu}

\author[M.~Bendersky]{M.~Bendersky}
\address{Department of Mathematics
CUNY,  East 695 Park Avenue New York, NY 10065, U.S.A.}
\email{mbenders@hunter.cuny.edu}

\author[F.~R.~Cohen]{F.~R.~Cohen}
\address{Department of Mathematics,
University of Rochester, Rochester, NY 14625, U.S.A.}
\email{cohf@math.rochester.edu}

\subjclass[2010]{Primary: 55P42, 55Q15, 52C35, 52B11, 35S22, 13F55;\\
Secondary:\/ 14M25, 14F45,  55U10, 55R20, 55N10}

\keywords{Polyhedral product, moment--angle complex, cohomology, arrangements, stable splitting, 
simplicial wedge, Davis--Januszkiewicz space, Golodness, monomial ideal ring}

\begin{abstract}
A polyhedral product is a natural subspace of a Cartesian product that is specified
by a simplicial complex. The modern formalism arose as a generalization of the spaces known as 
{\em moment--angle complexes\/} which were developed within the nascent subject of {\em toric topology\/}. 
This field, which began as a topological approach to toric geometry and aspects of symplectic geometry, has expanded 
rapidly in recent years. The investigation of polyhedral products  and their homotopy theoretic properties has 
developed to the point where they are studied in various fields of mathematics far removed from 
their origin.   In this survey, we provide a brief historical overview of the development of this subject, summarize many
of the main results and describe applications.
\end{abstract}
\maketitle
\tableofcontents
\section{Introduction}\label{sec:Introduction}
A polyhedral product is a natural topological subspace of a Cartesian product, determined by a simplicial complex
$K$ and a family of   pointed pairs of spaces  $(X_i,A_i)$, one for each vertex of $K$. As noted in \cite{th2},
in the special case that
$A_i$ is a basepoint for $X_i$, the polyhedral product mediates between  the product 
$X_1\times X_2\times\cdots \times X_m$ when $K$ is the full $(m-1)$-simplex, and the wedge 
$X_1\vee X_2 \vee \cdots \vee X_m$ when $K$ consist of $m$  discrete  points. 

 We begin with a brief historical note. The development of the theory of polyhedral products was guided by their 
inextricable link to spaces known as  moment--angle manifolds which arose within the subject of toric topology. This 
subject, which incorporates ideas from geometry, symplectic geometry, combinatorics and commutative algebra,
has now swelled well beyond its original confines as a topological generalization of 
toric geometry. 

Toric Topology has precipitated investigations into new areas of manifold and orbifold theory, centered 
around toric actions,  \cite{dj,bp1,bpr1,franz,bbko,gt1,hm,cray,prv,kmz,bss1,bnss}, and the unpublished notes of 
N.~Strickland,  \cite{strickland}. An excellent\ overview by V.~Buchstaber and N.~Ray is to be found in \cite{br}.

The study of polyhedral  products in particular has matured to the point that it extends now into a wide variety of
fields of mathematics distant from its origin, as detailed in the table below. Consequently, the development of their 
homotopy theory has become useful in this context. 

 In order to be more specific, it is now necessary to give a definition. Begin by setting $[m]   =\{1,2,\ldots,m\}$ and define 
a category $\bm{\mathcal{C}([m])}$  whose objects
are pairs $(\underline{X}, \underline{A})$ where
$$(\underline{X}, \underline{A}) \;=\; \{(X_1,A_1), (X_2,A_2), \ldots ,(X_m,A_m)\}$$ 
is a family of based CW-pairs. A morphism 
$$\underline{f}\colon (\underline{X}, \underline{A}) \longrightarrow (\underline{Y}, \underline{B})$$
consists of $m$ continuous maps $f_{i}\colon X_i \longrightarrow Y_{i}$ satisfying $f_i(A_i) \subset B_{i}$. 
Next, let $K$ be a simplicial complex on the vertex set $[m]$. We consider $K$ to be a category
where the objects are the simplices of $K$ and the morphisms $d_{\sigma,\tau}$ are the inclusions 
$\sigma \subset \tau$. 
\begin{defin}\label{defn:main}
 Let $K$ be a simplicial complex on the vertex set $[m]$ as above.  A  polyhedral product \index{Polyhedral product} is a functor 
$Z(K;-)\colon {\bm{\mathcal{C}([m])}} \longrightarrow {\rm {\bf Top}}$  satisfying
$$Z(K; (\underline{X}, \underline{A})) \; \subseteq \; \prod\limits_{i=1}^{m}{X_i}\;$$
\nd   that is given as the colimit of a diagram 
$D: K \to CW_{\ast}$, where at each $\sigma \in K$, we set 
\begin{equation}\label{eqn:d.sigma}
D(\sigma) =\prod^m_{i=1}W_i,\quad {\rm where}\quad
W_i=\left\{\begin{array}{lcl}
X_i  &{\rm if} & i\in \sigma\\
A_i &{\rm if} & i\in [m]-\sigma.
\end{array}\right.
\end{equation}
\end{defin}
\nd  Though here,  the colimit is a union given by 
$$Z(K; (\underline{X}, \underline{A})) = \bigcup_{\sigma \in K}D(\sigma),$$
 the full colimit structure is used heavily in the development of the elementary theory of polyhedral products. 
Notice that when $\sigma \subset \tau$  then $D(\sigma) \subseteq D(\tau)$.   In the case that $K$ itself is a simplex, 
$$Z(K; (\underline{X}, \underline{A}))\; =\; \prod\limits_{i=1}^{m}{X_i}.$$ 
\begin{rem}
When all the pairs $(X_i,A_i)$ are the same pair $(X,A)$, the family
$(\underline{X}, \underline{A})$ is written simply as $(X,A)$. Notice that in this case $Aut(K)$  acts naturally on 
 $Z(K;(X,A))$.
\end{rem}
\nd In a way entirely similar to that above, a related space $\widehat{Z}(K; (\underline{X}, \underline{A}))$, called the {\em smash
polyhedral product\/},  is defined by replacing the Cartesian product everywhere in Definition \ref{defn:main} by the 
smash product.
 That is,
$$\widehat{D}(\sigma) =\bigwedge ^m_{i=1}W_i \quad {\rm and} \quad
\widehat{Z}(K; (\underline{X}, \underline{A})) = \bigcup_{\sigma \in K}\widehat{D}(\sigma)$$
\nd  with
$$\widehat{Z}(K; (\underline{X}, \underline{A})) \; \subseteq \; \bigwedge\limits_{i=1}^{m}{X_i}.$$
 \index{Polyheral product, smash}

Polyhedral products generalize the spaces called {\em moment--angle complexes\/} \index{Moment--angle complex} which were
developed  by V.~Buchstaber and T.~Panov \cite{bp3} and correspond to the case
$$(X_{i},A_{i}) = (D^{2},S^{1}),  \quad i = 1,2,\ldots,m.$$  
 As mentioned above, the construction of moment--angle complexes  originated within the subject of  toric topology. 
Specifically, 
M.~Davis and T.~Januszkiewicz \cite{dj}, constructed spaces $\mathcal{Z}$ as  quotients of a product  of an 
$m$-dimensional torus
and an $n$-dimensional  
 polytope which is {\em simple\/}, meaning that it satisfies a certain homogeneity condition.
(See Section \ref{sec:genesis}.)  They were able to
realize a class of toric manifolds as quotients of $\mathcal{Z}$ by free torus actions.
The space $\mathcal{Z}$ was then reformulated as a moment--angle manifold by V.~Buchstaber and T.~Panov, \cite{bp3}. 

The moment--angle complex construction \eqref{eqn:d.sigmabp}  was extended to  pairs $(X,A)$ replacing
$(D^2,S^1)$ by  V.~Buchstaber and T.~Panov \cite {bp3}, N.~Strickland (unpublished), V.~Baskakov \cite{bask}, 
G.~Denham and A.~Suciu \cite{ds} and T.~Panov \cite{panxa}. The construction in the full generality of Definition 
\ref{defn:main} appeared in \cite{bbcg1}.
\newpage
The table below lists a selection of  current applications of polyhedral products.
\skp{0.3}
{\nd \hspace{0.3in}\begin{tabular}{ll}
$\bm{(\underline{X},\underline{A})}$    &$\bm{Z(K;(\underline{X},\underline{A})})$  \\
\\[-3mm]
$(D^{2},S^{1})$ &toric geometry and topology, arachnid mechanisms \\[0.3mm]
$(D^{1},S^{0})$ &surfaces, number theory, representation theory, linear programming\\[0.3mm]
$(S^{1}, \ast)$  &right--angled Artin groups,  robotics  \\[0.3mm]
$(\mathbb{R}\rm{P}^{\infty}, \ast)$ &right--angled Coxeter groups\\[0.3mm]
$(\mathbb{C}, \mathbb{C}^{\ast})$ &complements of coordinate arrangements \\[0.3mm]
$(\underline{\mathbb{R}}^{n}, (\underline{\mathbb{R}}^{n})^{\ast})$ &complements of certain non-coordinate arrangements\\[0.3mm]
$(\underline{\mathbb{C}\rm{P}}^{\infty}, \underline{\mathbb{C}\rm{P}}^{k})$\phantom{m} &monomial ideal rings\\[0.3mm]
$(\underline{EG}, \underline{G})$ &free groups and  monodromy  \\[0.3mm]
$(\underline{BG}, \underline{\ast})$ &combinatorics,  aspherical spaces  \\[0.3mm]
$(\underline{PX}, \underline{\Omega{X}})$ &homotopy theory, Whitehead products\\[0.3mm]
$(S^{2k+1}, \ast)$ &graph products, quadratic algebras
\end{tabular} 
\skp{0.4}
\nd\centerline{\LARGE{$\longrightarrow\!\longleftarrow$}}
\skp{0.1}
By way of motivation, we illustrate the Hopf map \index{Hopf map} and the Whitehead product \index{Whitehead product} 
from the point of view of polyhedral 
products. The $3$--sphere $S^3$ is realized  as a polyhedral product by taking $K$ to be
the simplicial complex consisting of two  discrete  points, 
$$\partial(D^2\times D^2) \; \cong\; S^1\times D^2 \cup_{S^1\times S^1}D^2 \times S^1 = Z(K;(D^2,S^1)) \;\subset\; 
D^2\times D^2 \subset \mathbb{C} \times \mathbb{C}.$$ 
The action of the compact two-torus $T^2$ on the product $D^2\times D^2$ leaves invariant the component pieces: 
$S^1\times D^2$, $S^1\times S^1$ and $D^2 \times S^1$, which comprise the basic building blocks of the 
polyhedral product, moment--angle complex, structure in this case.
The Hopf map 
$$S^3 \;\simeq\; S^1\times D^2 \cup_{S^1\times S^1} D^2 \times S^1 \longrightarrow D^2 \cup_{S^1} D^2\;\simeq\; S^2$$
is the quotient by the diagonal subgroup $S^1$ in $T^2$ acting in a natural way.
\skp{0.2}
 Consider next  the  Whitehead product $[f, g]$ defined by \index{Whitehead product}
$$S^3\; \stackrel{\omega}{\longrightarrow} \;S^{2}_\alpha\vee S^{2}_\beta \;   \xrightarrow{f\vee g} \; X$$
\nd  where $\omega$ is a map attaching the top cell of $S^2 \times S^2$ and $f\colon S^{2}_\alpha \longrightarrow X$
and $g\colon S^{2}_\beta \longrightarrow X$ are two basepoint preserving maps, \cite[page 381]{Hatcher:2001}.
The particular case $X = S^{2}_\alpha\vee S^{2}_\beta$ and $f = \iota_\alpha$, $g = \iota_\beta$, being the 
respective inclusions of $S^{2}_\alpha$ and  $S^{2}_\beta$  into $S^{2}_\alpha\vee S^{2}_\beta$, gives the Whitehead
product  $[\iota_\alpha, \iota_\beta]$. (Note here that $\iota_\alpha\vee \iota_\beta$ is the identity map.) 
From the point of view of polyhedral products\footnote{Lukas Katth\"{a}n alerted us to this fact.}, 
$[\iota_\alpha, \iota_\beta]$  is  induced by the map of pairs $(D^2, S^1)  
\longrightarrow (S^2, \ast)$  which collapses $S^1\subset D^2$ to a point,  
\begin{equation}\label{eqn:wp}
S^3 \;\simeq\; D^{2}_\alpha \times S^{1}_\beta \cup_{S^{1}_\alpha\times S^{1}_\beta} S^{1}_\alpha\times D^{2}_\beta  
\longrightarrow  (S^{2}_\alpha \times \ast) \cup_{\ast \times \ast} (\ast \times S^{2}_\beta) \;\simeq\; S^{2}_\alpha \vee S^{2}_\beta.
\end{equation}
Here, the space on the right-hand side is the polyhedral product $Z(K;(\underline{S}^2,\ast))$ where again, $K$ consists
of two  discrete  points. 
Within this subject, wedge products appear often, 
thus Whitehead products appear in profusion. These topics are discussed in more detail in Section \ref{sec:hwp}.
\skp{0.1}
\nd  \centerline{\LARGE{$\longrightarrow\!\longleftarrow$}} 
\skp{0.1}
 The paper is organized as follows. A brief technical discussion  of the origin of polyhedral products as moment--angle 
complexes in  
Sections \ref{sec:genesis} and \ref{sec;bpapproach}, includes the original topological construction of toric manifolds 
by M.~Davis and
T.~Januszkiewicz \cite{dj} and the reformulation of one of their main spaces as a moment--angle complex by V.~Buchstaber and 
T.~Panov, \cite{bp3, bp1}.

A  discussion follows in Section \ref{sec:santiago} about the independent discovery of moment--angle complexes
as intersections of certain quadrics. We describe briefly the work of  S.~L\'opez de Medrano \cite{ldm1},  
S.~L\'opez de Medrano and A.~Verjovsky \cite{ldmv}, F.~Bosio and L.~Meersseman \cite{bm},  followed by that of 
 S.~L\'opez de Medrano  and S.~Gitler \cite{gl}. This approach has proved effective in identifying classes of moment--angle 
 manifolds that are connected sums of   products of spheres.   This is followed in Section \ref{sec:cohomofamac} by a 
 sketch of the computation of the cohomology ring of moment--angle complexes.

We begin our focus on more general polyhedral products in Section \ref{sec:exp}, by describing their behaviour
with respect to the exponentiation of CW-pairs. This has led to  an application to the study of toric manifolds 
\cite{erokhovets,u1,u2,suciu,gt2,gptw,cp1},  and orbifolds \cite{bss2}. An application to topological joins is included as 
an example of the utility of this property with respect to exponentiation.

The behaviour of polyhedral products with respect to fibrations, due to G.~Denham and A.~Suciu \cite{ds}, is surveyed 
briefly in Section \ref{sec:early}.
We take advantage of the formalism to describe briefly the early work in this area by G.~Porter \cite{porter, porter2},
T.~Ganea \cite{ganea} and A.~Kurosh \cite{kurosh}. and G.~W.~Whitehead.  Instances of polyhedral products 
appear also in the work of D.~Anick \cite{anick}.

In Section \ref{sec:splitting}, we review the various fundamental unstable and stable splitting theorems for the polyhedral
product. These theorems give access to the homotopy type of the polyhedral product in many of the most important
cases and drive a number of applications.  The identification of the various wedge summands which appear in the stable
decompositions, occupies the second part of this section.

A fine theorem of A.~Al-Raisi \cite{ar} shows that,  in certain cases,  the stable splitting of the 
polyhedral product can be chosen to be 
equivariant with respect to the action of $Aut(K)$. This and other related results, described in Section \ref{sec:equivariance},
represent the first foray of the subject into representation theory. An  application of Al-Raisi's theorem to surfaces yields 
a cyclotomic identity.

The important case \apoint\; is discussed in Section \ref{sec:aisapoint}. It arises in algebraic combinatorics, the study of free
groups and monodromy representations, geometric group theory, right--angled Coxeter and Artin groups and asphericity,

In this section also, we use a recent result of T.~Panov and S.~Theriault \cite{pt} to give a short 
proof that if $K$ is a flag simplicial complex then, for a discrete group $G$,  $Z(K;(BG,\ast))$ is an Eilenberg--Mac Lane 
space, a known result, \cite{chardav,do,ms2}. 

Results on the cohomology of polyhedral products are presented in Section \ref{sec:cohomology}. The discussion begins
with a description of the cup-product structure as a consequence  of the stable splitting 
Theorem \ref{thm:bbcgsplitting}. A spectral sequence corresponding to a natural filtration of $Z(K;\xa)$ by the simplices of
$K$ is described next and followed by several examples.

Section \ref{sec:cartan} describes an entirely geometric approach to the calculation of the cohomology. The notion of 
{\em wedge decomposable pairs\/} is introduced and the corresponding polyhedral products are determined  explicitly
as wedges of identifiable spaces. Moreover, it is shown in Theorem \ref{thm:cartan2}, that over a field,  every CW-pair $\xa$ 
of finite type is cohomologically wedge decomposable. This  suffices to give the additive structure of the 
homology of any $Z(K;(X,A))$ for CW-pairs $(X,A)$ of finite type, and any field coefficients.

The application of polyhedral products to questions concerning the Golod properties of certain rings is 
described in 
Section \ref{sec:golodness}. Higher Whitehead products, constructed using polyhedral products, are discussed in
Section \ref{sec:hwp}. In cases when a moment--angle complex splits into a wedge of spheres, higher  Whitehead
products describe certain canonical maps related to the inclusions of the summands, \cite{pr,gt3, ik6,abramyan}. 
They arise also in the analysis of $\Omega{Z}(K;(D^2,S^1))$.

Though the $KO$--theory of certain polyhedral products and toric manifolds is discussed in the literature 
\cite{bbko, kodj}, it is omitted from this article. Applications to robotics (\cite{hck, kt}, for example)  are also omitted.

\section{The origin of polyhedral products in toric topology}\label{sec:genesis} 
 A {\em polytope\/} is the convex hull of a finite set of points in some $\mathbb{R}^n$. The
{\em  dimension\/}  of the polytope is the dimension of its affine hull, \cite[Chapter 1]{bp1}. 
We shall assume that an $n$-dimensional polytope is a subset of $\mathbb{R}^n$. A codimension-one face of a polytope
is called a {\em facet\/}. A polytope $P^n$ is called {\em simple\/} if at each vertex, exactly $n$ facets intersect. The
convention is to let $m$ denote the number of facets of $P^{n}$.

In order to introduce the notion of  a {\em toric manifold\/}, we recall that the real  torus $T^n$ acts
on $\mathbb{C}^{n}$ in a standard way, and the quotient is 
$$\mathbb{C}^{n}\big/T^{n} \;\cong\; 
\mathbb{R}^n_+ \;=\; \{ (x_1,\ldots, x_n) \in \mathbb{R}^n: x_i \geq 0 \mbox{ for } i=1, \ldots , n\}.$$ 

 A toric manifold $M^{2n}$, (sometimes called a {\em quasitoric manifold\/} in the literature), is a compact manifold 
covered by local charts $\mathbb{C}^n$, with a global action of a real $n$-dimensional torus
$T^n$, that restricts to the standard action on each $\mathbb{C}^n$ component. Under this action, each copy of 
$\mathbb{C}^n$ must have a quotient $\mathbb{C}^{n}\big/T^{n} \;\cong\; \mathbb{R}^n_+ $ which is a neighbourhood
of a vertex of a simple polytope $P^n \cong M^{2n}\big/T^n$. 
\index{Toric manifold}
\begin{rem}
Smooth projective toric varieties are examples of toric manifolds in the sense above.
\end{rem}
A topological approach
to the construction of toric manifolds was developed by M.~Davis and T.~Januszkiewicz in \cite{dj}. 
They begin with a simple polytope $P^n$ and a function
$$\lambda\colon {\mathcal F} = \{F_{1},F_{2},\ldots,F_{m}\} \longrightarrow \mathbb{Z}^n.$$
from the set of facets of $P^n$ 
\nd into an $n$-dimensional integer lattice. The function $\lambda$ must satisfy a {\em regularity\/}
condition, namely, whenever 
$F = F_{i_{1}} \cap F_{i_{2}} \cap \cdots \cap F_{i_{l}}$ is a codimension-$l$ face of $P^n$, then 
$\{\lambda(F_{i_{1}}),\lambda(F_{i_{2}}), \ldots ,\lambda(F_{i_{l}})\}$ spans 
an $l$-dimensional submodule of $\mathbb{Z}^n$ which is a direct summand. (The fact that every such
face can be written uniquely as such an intersection is a consequence of the polytope being simple.) 
Such a map is called a {\em characteristic function\/} \index{Characteristic function} associated to $P^n$.

Next, regarding $\mathbb{R}^n$ as the Lie algebra of $T^n$,  the map
$\lambda$ is used to associate  to each codimension-$l$ face $F$ of $P^n$ \/ a rank-$l$
subgroup $G^\lambda_F \subset T^n$, namely the subgroup of $T^n$ determined by the span
of the $\lambda(F_{i_{j}})$. (Specifically, writing $\;\lambda(F_{i_j}) = (\lambda_{1{i_j}},\lambda_{2{i_j}},\ldots,\lambda_{n{i_j}})$
\nd gives
\begin{equation}\label{eqn:gf}
G^\lambda_F = \big\{\big(e^{2\pi{i}(\lambda_{1{i_1}}t_1 +\lambda_{1{i_2}}t_2 + \cdots + \lambda_{1{i_l}}t_l)},
\ldots, e^{2\pi{i}(\lambda_{n{i_1}}t_1 + \lambda_{n{i_2}}t_2 + \cdots +\lambda_{n{i_l}}t_l)}\big) \in T^n\big\}
\end{equation}
\nd where $t_i \in \mathbb{R},\, i = 1,2,\ldots,l$.) 

Finally, let $p \in$ $P^n$  and $F(p)$ be 
the unique face with $p$ in its relative interior. Define an equivalence
relation $\sim_{\!\lambda}$ on $T^n$ $\times$ $P^n$ $\;$ by $(g,p) \sim_{\!\lambda} (h,q)$ if and only
if $p = q$ and $g^{-1}h \in G^\lambda_{F(p)} \cong T^l$. Then
\begin{equation}\label{eqn:defn.tn}
M^{2n} \cong M^{2n}(\lambda) = T^n \times P^n\big/\!\sim_{\!\lambda}
\end{equation} 
\nd is a smooth, closed, connected, $2n$-dimensional toric manifold with 
$T^n$ action induced by left translation \cite[page 423]{dj}. A projection 
$$\pi \colon M^{2n} \rightarrow P^n \cong M^{2n}\big/T^n$$
\nd onto the polytope is induced from the projection
$T^n \times$ $P^n$ $\rightarrow$ $P^n$. 
\begin{exm}
Consider  the case for which  $P^n = [0,1]$, the one-simplex having dimension $n=1$ and two facets, so $m=2$. 
Here, the set of facets ${\mathcal F} = \{F_{1},F_{2}\}$ consists of the two vertices of the one-simplex.
Define $\lambda : {\mathcal F} \longrightarrow \mathbb{Z}^1$ by $\lambda(F_1) = -1$ and $\lambda(F_2) = 1$.
Topologically, $T^1 \times P^1$
is the cylinder $T^1 \times [0,1]$. In this case, the equivalence relation $\sim_{\!\lambda}$ collapses to points 
each circle $T^1 \times \{0\}$ and $T^1 \times \{1\}$, yielding
\begin{equation*}
T^1 \times P^1\big/\!\sim_{\!\lambda}\;\;\cong\; S^2.
\end{equation*} 
\nd In fact, this $S^2$ has the structure of the toric variety $\mathbb{C}\rm{P}^1$.
\end{exm}
\begin{rem}
Every projective non-singular toric 
variety has the description \eqref{eqn:defn.tn}; $P^n$ \;and $\lambda$ encode topologically the information in the defining fan, 
\cite[Chapter 5]{bp2}.
\end{rem}

As part of their approach to the construction of toric manifolds, the authors of \cite{dj} 
introduced  a {\em second\/} manifold  
\begin{equation}\label{eqn:mac}
\mathcal{Z}  = T^m \times P^n \big/\!\!\sim
\end{equation}
\nd where $T^m$ denotes the real $m$--torus, its coordinate circles indexed by the $m$ facets of the simple polytope
$P^n$.  (Notice here a distinction between \eqref{eqn:mac} and \eqref{eqn:defn.tn} given by the fact that
$T^m$ has replaced $T^n$; the equivalence relation will be different too. This space was introduced into the 
theory of toric manifolds in \cite{dj} in order to facilitate the calculation of the cohomology of the toric manifold
$M^{2n}(\lambda)$ in \eqref{eqn:defn.tn}.)  The equivalence relation $\sim$ is given by defining 
$\theta\colon {\mathcal F} \longrightarrow \mathbb{Z}^m$ by $\theta(F_{i}) = \underline{e}_{i} \in \mathbb{Z}^{m}$, 
the standard basis vector. 
\begin{rem}
 Unlike a characteristic map $\lambda$ (above) the map $\theta$ indexes the coordinate circles in $T^m$ by the 
facets of $P^n$, {\em and so depends on the combinatorics of the polytope $P^n$ only\/}.
\end{rem}
Next, by exact analogy with the construction of the subgroups $G^{\lambda}_{F} \subset T^n$ above for each
codimension-$l$ face $F$ of $P^n$, the map $\theta$ determines a rank-$l$ subgroup $G^{\theta}_{F} \subset T^m$.
The explicit description of this group is similar to \eqref{eqn:gf}. 
Finally, the equivalence relation $\sim$ is defined in a manner entirely analogous to $\sim_{\lambda}$
above. That is, for  $p \in$ $P^n$, and $F(p)$, the unique face with $p$ in its relative interior,  define 
$\sim$ on $T^n$ $\times$ $P^n$ $\;$ by $(g,p) \sim (h,q)$ if and only
if $p = q$ and $g^{-1}h \in G^\theta_{F(p)} \cong T^l$. 

It follows that $\mathcal{Z}$ is a smooth, closed, connected, $(m+n)$--dimensional manifold endowed with a
$T^m$ action induced by left translation \cite[page 423]{dj}. A 
projection $\mathcal{Z}  \rightarrow \mathcal{Z} \big/T^m \cong P^n$ onto the polytope is induced from the projection
$T^m \times$ $P^n$ $\rightarrow$ $P^n$.  (See also \cite{atiyah} and \cite{jurk}.) 
\begin{exm}
Consider the case for which $P^n$ is the one-simplex again so that $n=1$ and $m=2$. Topologically, 
$T^m \times P^n$
is the cylinder $T^2 \times [0,1]$. The equivalence relation collapses one of the coordinate circles
in $T^2 \times \{0\}$  to a point as a group, as well as the {\em other\/} coordinate circle in $T^2 \times \{1\}$. This gives 
$$\mathcal{Z} \;\cong\; S^1 \ast S^1 \;\cong\; S^3.$$
\end{exm}


A reformulation of the ideas above indicates that every toric manifold $M^{2n}$ can be described as a quotient
\begin{equation}\label{eqn:zquotient}
M^{2n} \cong M^{2n}(\lambda)  \; \cong \; \mathcal{Z}\big/T^{m-n}
\end{equation}
\nd where here,  $T^{m-n} \subset T^m$ is the torus determined by the kernel of the map $\lambda$, which the regularity
condition  of $\lambda$  ensures acts freely.  In the next section, the space $\mathcal{Z}$ is reformulated as a 
{\em moment--angle manifold\/}.

\section{The introduction of moment--angle complexes}\label{sec;bpapproach}
Every simple polytope $P^n$ with $m$ facets, embeds naturally into   the union of  the facets of the 
$m$-cube $I^m$ by virtue of its
{\em cubical decomposition\/}, \cite[Construction $4.5$]{bp1}. V.~Buchstaber and T.~Panov constructed a $T^m$--invariant 
natural subspace $Z_{K_{P^n}}$ of 
$\mathbb{C}^m$, 
that completes a commutative diagram:
\begin{equation}\label{eqn:cubical}
\begin{CD}
Z_{K_{P^n}} @>{\alpha}>>
(D^2)^m\subset \mathbb{C}^m\\
@VV{/T^m}V                                          @VV{/T^m}V \\
P^n
@>{\iota}>>
I^m
\end{CD}
\end{equation}
\nd Here, $D^2 \subset  \mathbb{C}$   is the unit disc with the standard circle action, the vertical maps are
quotients by the standard action of $T^m$,  and $K_{P^n}$ denotes the boundary
complex of the dual of $P$, 
 (the vertices of $K_{P^n}$ are the facets of $P^n$ and simplices correspond to intersections of facets),  
\cite[Example 2.2.4]{bp2}. They recognized the space $Z_{K_{P^n}}$ as the natural subspace of 
$(D^2)^m$  defined as a colimit by a diagram 
$D: K_{P^n} \to CW_{\ast}$, which at each $\sigma \in K_{P^n}$, 
is given by
\begin{equation}\label{eqn:d.sigmabp}
D(\sigma) =\prod^m_{i=1}W_i,\quad {\rm where}\quad
W_i=\left\{\begin{array}{lcl}
D^2 &{\rm if} & i\in \sigma\\
S^1 &{\rm if} & i\in [m]-\sigma.
\end{array}\right.
\end{equation}
\nd Here, the colimit is $Z_{K_{P^n}} \;=\; \bigcup_{\sigma \in K_{P^n}}{D(\sigma)}$, and is now written
$Z(K_{P^n};(D^2,S^1))$; it was called a {\em moment--angle manifold\/} by V.~Buchstaber and T.~Panov. 
\index{moment--angle manifold}

In order to  verify that the moment--angle complex $Z(K_{P^n};(D^2,S^1))$ is in  fact the space $\mathcal{Z}$ constructed
by M.~Davis and T.~Januszkiewicz, we work locally.
Construction \eqref{eqn:mac} can be analyzed in the neighborhood of a vertex, where the simple polytope 
$P^{n}$ is PL-homeomorphic to $\mathbb{R}^{n}_{+}$.
\newpage
\nd \hspace{1.7in}\includegraphics[width=2.8in]{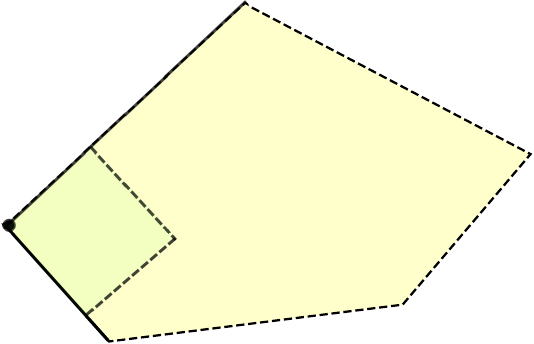}
\skp{-2.05}\hspace{1.3in}{\large $v_i$}\hspace{0.37in}{\large $\mathbb{R}^{n}_{+}$}
\skp{-0.3}\hspace{3.1in}{\large $P^n$}
\skp{1.5}
\centerline{{\sc Figure 1.} The local structure of a simple polytope $P^n$}
\skp{0.2}
\nd The polytope can be given a {\em cubical\/} structure as in 
\cite[Construction $5.8$ and Lemma $6.6$]{bp1}. A cube $I^{n}$, {\em anchored\/} by the vertex $v_{i}$, sits
inside the copy of $\mathbb{R}^{n}_{+}$  which arises  by deleting all faces of $P^{n}$ which do not contain $v_{i}$.
Locally, $T^{m}\times P^{n}$ is
\begin{equation}\label{eqn:local}
T^{m}\times I^{n} \;\cong\; (S^{1} \times I)^{n} \times (S^{1})^{m-n}.
\end{equation}
\nd Recall that the map $\theta$, defined in Section \ref{sec:genesis}, indexes all the
circles in $T^{m}$ by the facets of the polytope, so the order of factors here has 
been shuffled naturally. The factors $S^{1}$ which are paired  with a copy of $I$ are those 
corresponding to the facets of $P^{n}$ which meet at $v_{i}$. The effect of the equivalence relation $\sim$  in \eqref{eqn:mac}
on $T^{m}\times I^{n}$
is to convert every $S^{1} \times I$ on the right-hand side of \eqref{eqn:local} into a 
disc $D^2$, by collapsing $S^{1}\times \{0\}$ to a point. So
\begin{equation}\label{eqn:blocks}
T^{m}\times I^n\big/\!\!\sim\quad\cong\quad (D^{2})^{n} \times (S^{1})^{m-n}.
\end{equation}
\nd The vertices of $P^{n}$ correspond to the maximal simplices of the simplicial complex $K_{P^n}$ so, assembling the blocks 
\eqref{eqn:blocks}, gives the moment-angle manifold
\begin{equation}\label{eqn:zequalszk}
\mathcal{Z} \;=\; T^{m}\times P^{n}\big/\!\!\sim \quad\cong\quad Z(K_{P^n}; (D^{2},S^{1})).
\end{equation}
\begin{remark}\label{rem:generalk}
In a natural way, the Buchstaber--Panov formulation generalizes from $K_{P^n}$ to any simplicial complex $K$, but the result 
is not in general a manifold. For general $K$ the space $Z(K;(D^2,S^1))$ is called a moment--angle {\em complex\/}. 
The condition that $K$ be a polytopal sphere, \index{Polytopal sphere}
  (a triangulated sphere which is isomorphic to the boundary complex of  
a polytope which has all its faces simplices, a simplicial polytope, \cite[Definition $2.5.7$]{bp2}),  
ensures that $Z(K;(D^2,S^1))$ is a manifold. This was weakened in 
 \cite[Corollary $2.10$]{cai1} 
to: the moment--angle complex is a topological $(n +m)$--manifold if and only if the realization of $K$ is a 
{\em generalized homology 
$(n-1)$--sphere\/}.  The latter is defined to be a triangulated polyhedron, which has the homology of an $(n-1)$--sphere for
which the link of each $p$-simplex has the homology of an $(n-2-p)$--sphere \cite{cai1}. \index{Generalized homology sphere}

 Currently, the best result implying that the moment--angle complex is a smooth manifold, requires that  $K$ be the simplicial 
complex underlying a complete simplicial fan;  more details are to be found in 
\cite[Theorem $2.2$]{pu} and \cite[Theorem $9.2$]{panovman}.

 The construction of moment--angle complexes was extended to simplicial posets by Z.~L\"{u} and T.~Panov in
\cite{lupan} where they are used to study the face rings of posets topologically. 
\end{remark}
Below, the $2$-torus is exhibited as the polyhedral product $Z(K;(D^1,S^0)) \subset (D^1)^4$, where $K$ is the boundary 
complex of a square on vertices
$\{1,2,3,4\}$. According to Definition \ref{defn:main}, the empty simplex of $K$ contributes $(S^0)^4$ which is identified with
the $16$ vertices of $(D^1)^4$, and the four vertices of $K$ contribute the
one-skeleton consisting of 32 edges. The maximal simplices of $K$ are: $\sigma_1 = \{12\}, \sigma_2 = \{14\}, \sigma_3 = \{23\},
\sigma_4 = \{34\}$. 
\begin{rem}
For clarity, we have included  the whole one-skeleton in each $D(\sigma_i)$ below, but, strictly speaking,
the {\em whole\/} one-skeleton should appear in the final picture only.
\end{rem} 
\skp{0.3}
\nd \hspace{0.2in}\includegraphics[width=1in]{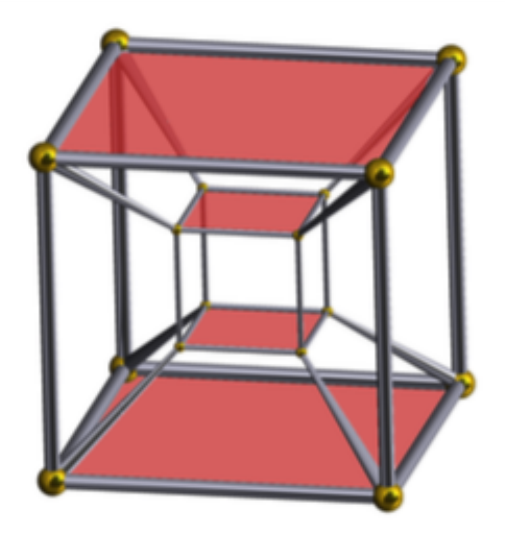} \hspace{0.1in}\includegraphics[width=1in]{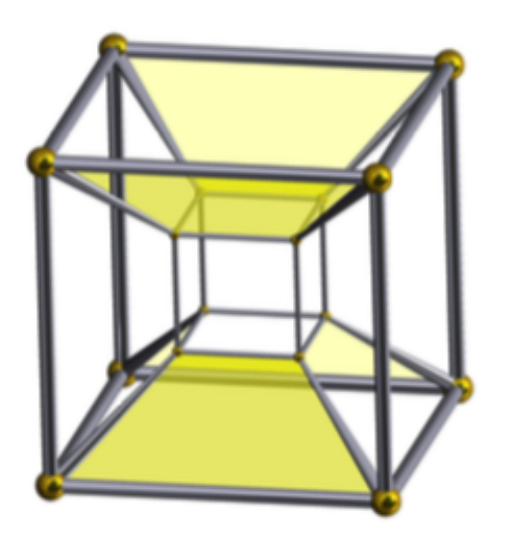}
\nd \hspace{0.2in}\includegraphics[width=1in]{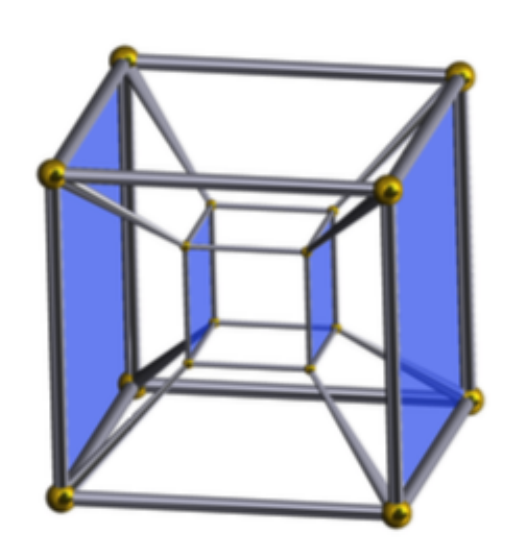} \nd \hspace{0.2in}\includegraphics[width=1in]{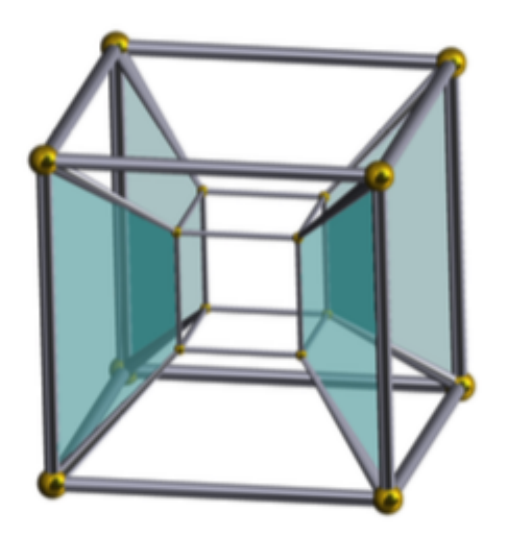}
\nd \hspace{0.2in}\includegraphics[width=1in]{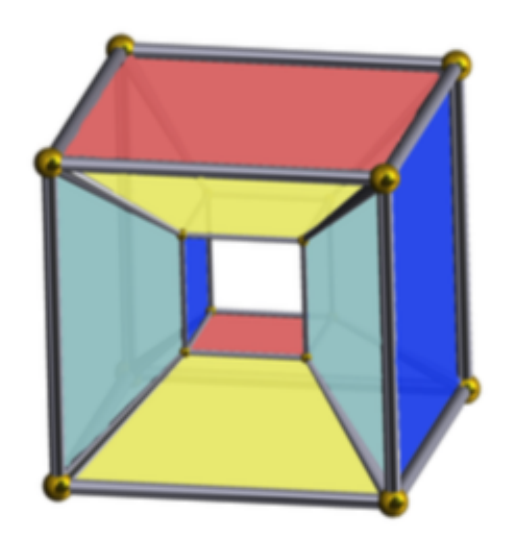}
\skp{0.02} 
\nd \hspace{0.4in} $D(\sigma_{1})$ \hspace{0.6in}  $D(\sigma_{2})$ \hspace{0.7in} $D(\sigma_{3})$ \hspace{0.7in} 
$D(\sigma_{4})$ \hspace{0.7in} the torus 
\skp{0.1}
\centerline{{\sc Figure 2.} The $2$--torus as the polyhedral product $Z(K;(D^1,S^0))$}
\skp{0.01}
\centerline{(This figure is reproduced by kind permission of Alvise Trevisan, cf.~\cite{at}.)}

\section{Moment-angle complexes as intersections of quadrics}\label{sec:santiago}  \index{Quadrics, intersection of}
The purpose of this section is to provide an overview of the independent topological development of moment--angle complexes
arising from beautiful work of F.~Bosio and L.~Meersseman \cite{bm}, Y.~Barreto, S.~L\'opez de Medrano and A.~Verjovsky 
\cite{bldmv},  S.~L\'opez de Medrano \cite{ldm1},  as well as his subsequent joint work with S.~Gitler \cite{gl},  and
S.~L\'opez de Medrano and   A.~Verjovsky \cite{ldmv},   Closely related work was developed earlier, 
and independently, by C.~T.~C.~Wall \cite{wall}. All of these 
authors considered quadrics specified by certain natural choices of homogeneous quadratic polynomials intersected with 
the unit sphere. One consequence is the L\'opez de~Medrano--Gitler proof of a conjecture of Bosio--Meersseman  in the case that
the polytope $P^n$ is even dimensional and  ``dual neighborly''.  For such polytopes $P^n$, they show that 
$Z(K_{P^n};(D^1,S^0))$ is a connected sum of  products of spheres.   In addition to the proof of this conjecture, 
the authors of \cite{gl} 
also give further computations of various cup-products in moment-angle manifolds.

 Let $m$ and $n$ be positive integers satisfying $m>n$.  Set $k=2m-n-1$. A  certain  manifold, 
 $M \subset {\mathbb{C}}^m$ 
 associated to an $n$-dimensional polytope, is defined  next  as the intersection of $k$ quadrics.  Specifically with $\lambda_i \in 
 \mathbb{R}^k,  i= 1, \ldots ,m$\; the manifold $M$ is the intersection of the solution of a  homogeneous  system of
 equations with the unit sphere in $\mathbb{C}^m$,
$$ \Big\{(z_1,z_2,\ldots,z_m) \in \mathbb{C}^m\colon\underset{i=1}{\overset{m}{\sum}} \lambda_i |z_i|^2=0\Big\}\;\;
\cap \;\;\Big\{(z_1,z_2,\ldots,z_m) \in \mathbb{C}^m\colon \underset{i=1}{\overset{m}{\sum}} |z_i|^2=1\Big\}.$$ 
The manifold M admits an action of the $m$--torus $T^m$ by multiplication
of unit complex numbers on each of the coordinates.  The quotient by this action can be identified with the $n$-dimensional 
convex  polytope $P^n$ given by
$$\underset{i=1}{\overset{2m}{\sum}} \lambda_i x_i =0,\quad \underset{i=1}{\overset{2m}{\sum}} x_i =1,\quad x_i \geq0.$$
This construction is essentially equivalent to the construction \eqref{eqn:mac} of a moment-angle manifold from a 
simple polytope in Section \ref{sec:genesis}. The simplicity of $P^n$ is implied by  the weak hyperbolicity condition: 
\index{Hyperbolicity, weak} {\it the convex hull of any subset of $\{\lambda_1, \ldots , \lambda_m\}$, with
 $k$ or fewer elements does not contain the origin}.
The real points of this manifold correspond to the associated real moment-angle manifold $Z(K_{P^n};(D^1,S^0))$ that has 
natural a $\mathbb{Z}_{2}^n$-action.

S.~L\'opez de Medrano and A.~Verjovsky \cite{ldmv}  constructed  a class of non-K\"ahler, compact, manifolds generalizing 
Hopf and Calabi-Eckmann manifolds that are projectivizations of some of these moment--angle manifolds. L.~Meersseman 
\cite{meers}, generalized this construction and also observed that all even-dimensional moment--angle manifolds are 
themselves in this class. 

 S.~L\'opez de Medrano \cite{ldm1}, showed that when $k=2$,  the moment--angle manifolds constructed above and their 
 real parts  are either three-fold products  of spheres or connected sums of two-fold products of spheres.

These works culminated in the article by S.~Gitler and S.~L\'opez de Medrano \cite{gl}, showing that for all values of $k$,  these 
manifolds are  frequently connected sums of products of spheres. 

One result is that if $K$ is a polytopal sphere,   (see Remark \ref{rem:generalk}) 
then any manifold $Z(K;(D^1,S^0))$ which is  $2c$-dimensional, 
and  $(c-1)$-connected for $c\geq 3$, must be diffeomorphic to a connected sum of 
 products of spheres.   
This is a proof of a conjecture of Bosio and  Meersseman  \cite[page 115]{bm}, in the case of
even dimensional dual-neighborly polytopes mentioned earlier in this section. 
The authors also introduce operations including {\em corner cutting\/} or {\em edge cutting\/} of $K$ and  they analyze the effect of 
these operations on the manifolds $Z(K;(D^1,S^0))$.  One result in this direction is that they obtain new, infinite families of 
manifolds  which are diffeomorphic to connected sums of products of spheres.

The methods in \cite{gl} influence the argument used by S.~Theriault in \cite{th1} to prove Panov's Conjecture. 
\index{Panov's conjecture} The
number of wedge summands in a homotopy decomposition of the moment-angle complex   corresponding to the simplicial
complex consisting of $l$  discrete  points, is described in \eqref{eqn:disjointpoints}. Panov's conjecture relates this 
to the connected sum factors in a decomposition of the moment--angle manifold,  up to diffeomorphism, 
corresponding to the simple polytope obtained by making $l$  corner cuts   on a standard simplex.

 The varieties   discussed above, appearing in work of S.~L\'opez de Medrano et alia, split after suspending once. 
Other varieties, given by classical Stiefel manifolds, also stably decompose by work of I.~M. James \cite{james2}, 
and H.~Miller \cite{haynes}. It is 
natural to consider the variety given by intersections of these two types of quadratic varieties, as well as examine whether, 
and how they stably decompose. A result of Cohen and Peterson \cite{copet} shows that the stable decompositions for 
Stiefel varieties generically require many suspensions in order to split. 
\section{The cohomology of moment--angle complexes}\label{sec:cohomofamac}
The integral cohomology ring of a moment--angle complex was computed  by  M.~Franz \cite{franz}
and independently by V.~Baskakov, V.~Buchstaber and T.~Panov \cite{bbp}. 
In the latter, the two-disc $D^2\subset \mathbb{C}$ is decomposed into three cells: the zero-cell is the point
$1\in D^2$, the one-cell, denoted by $T$, is the complement $\partial{D^2} \smallsetminus\{1\}$ and the
two-cell, denoted by $D$ is the interior of $D^2\subset \mathbb{C}$. Taking products yields a cellular decomposition
of $D^m\subset \mathbb{C}^m$. 

Below we summarize the computation, following the development in \cite{bp2}.  The cells of $D^m$ are parametrized 
by pairs of subset $I, J \subset  [m] = \{1,2,3,\ldots\}$  satisfying $I\cap J = \varnothing$. 
The set $J$ parametrizes the $T$-cells and the set 
$I$ parametrizes the $D$-cells. The remaining  positions $[m] \smallsetminus(I\cup J)$ are occupied by $0$-cells. The cell 
corresponding to the pair $J,I$ is denoted $\langle J,I\rangle$. Next, we introduce $K$, a simplicial complex 
 on the vertices $[m]$. 
In this decomposition, the moment--angle complex $Z(K;(D^2,S^1))$ embeds in $D^m$ as a cellular subcomplex;
the cell $\langle I,J\rangle$ is in $Z(K;(D^2,S^1))$ whenever $I$ corresponds to a simplex in $K$.

The cellular cochains $C^{\ast}\big(Z(K;(D^2,S^1))\big)$ have a basis $\langle I,J\rangle^\ast$ dual to the cells. These
are bigraded by deg$\langle I,J\rangle^\ast = (-|J|, \;2|I| + 2|J|)$. The cellular differential preserves the {\em second\/} grading
and so the complex splits as follows:
$$C^{\ast}\big(Z(K;(D^2,S^1))\big) = \bigoplus_{q=0}^{m}C^{\ast,2q}\big(Z(K;(D^2,S^1))\big).$$
The {\em bigraded Betti numbers\/} are defined by 
$$b^{-p,2q}(K) \;=\; {\rm rank}\hspace{0.6mm}H^{-p,2q}\big(Z(K;(D^2,S^1))\big)$$
and so the ordinary Betti numbers of the moment-angle complex satisfy
$$b^k(K) = \sum_{-p+2q=k}b^{-p,2q}(K).$$
 
Recall that the integral Stanley--Reisner ring of $K$, $\mathbb{Z}(K)$, is a polynomial ring on two-dimensional 
generators $v_i$, modulo the ideal  $I$ generated by monomials corresponding to the minimal non-faces of $K$,
\index{Stanley--Reisner ring} \index{Face ring}
\begin{equation}\label{eqn:srring}
\mathbb{Z}(K) \cong \mathbb{Z}[v_1,\ldots,v_m]/I.
\end{equation}
 A {\em minimal non-face\/} of 
$K$ is a sequence of vertices of $K$ which is not a simplex of $K$, but any proper subset is a simplex of $K$. 

\nd There is an isomorphism of cochain complexes, \cite[Lemma $4.5.1$]{bp2},
$$g\colon R^{\ast}(K) \longrightarrow C^{\ast}\big(Z(K;(D^2,S^1))\big)$$
\nd where $R^{\ast}(K) $ is the quotient of a Koszul algebra
$$\Lambda[u_1,\ldots,u_m] \otimes \mathbb{Z}(K)\big/(v_{i}^2 = u_{i}v_i = 0, 1\leq i\leq m)$$
\nd with bigrading and differential given by 
\begin{equation}\label{eqn:rofk} 
{\rm deg}\/ u_i = (-1,2), \;\; {\rm deg} v_i = (0,2),\;\;du_i = v_i,\;\; dv_i=0.
\end{equation}
\nd The isomorphism $g$ is given by $u_{J}v_I \mapsto \langle I,J\rangle^\ast$. It follows that there is an additive 
isomorphism
$$H\big(R^{\ast}(K)\big) \cong H^\ast\big(Z(K;(D^2,S^1))\big).$$

The properties of polyhedral products are used then to show that the algebra $R^{\ast}(K)$ is weakly equivalent to
the Koszul algebra as follows;
$$\Lambda[u_1,\ldots,u_m] \otimes \mathbb{Z}(K) = C^{\ast}\big(Z(K;(S^\infty,S^1))\big) \cong 
C^{\ast}\big(Z(K;(D^2,S^1))\big) = R^{\ast}(K).$$
Furthermore, a  cellular diagonal approximation is constructed in \cite{bp2}, which allows for an extension to the 
multiplicative structure of the cohomology. \index{Moment--angle complex, cohomology}
\begin{thm} \cite[Theorem $4.5.4$]{bp2}\label{thm:cohomofz}
There are isomorphisms of bigraded algebras
$$H^{\ast,\ast}\big(Z(K;(D^2,S^1))\big) \cong {\rm Tor}_{\mathbb{Z}[v_1,\ldots,v_m]}\big(\mathbb{Z}(K),\mathbb{Z}\big)
\cong H\big(\Lambda[u_1,\ldots,u_m] \otimes \mathbb{Z}(K) , d\big)$$
where the bigrading and differential are as in \eqref{eqn:rofk}. The isomorphisms are natural for simplicial 
morphisms of $K$,  which are assumed to be monomorphisms.
\end{thm}
\nd Notice now that the bigraded Betti numbers satisfy \index{Betti numbers, bigraded}
$$b^{-p,2q}(K) = {\rm dim\hspace{0.5mm}Tor}^{-p,2q}_{\mathbb{Z}[v_1,\ldots,v_m]}\big(\mathbb{Z}(K), \mathbb{Z}\big)$$
\nd and their calculation is simplified by the following theorem due to M.~Hochster. \index{Hochster's Theorem}
\begin{thm}[\cite{hochster}]\label{thm:hochster}
Let $K_J$ denote the full subcomplex of $K$ corresponding to $J \subset [m]$. Then
$${\rm Tor}^{-p,2q}_{\mathbb{Z}[v_1,\ldots,v_m]}\big(\mathbb{Z}(K), \mathbb{Z}\big)\;\;= 
\bigoplus_{J\subset [m],  |J| =q}\widetilde{H}^{q-p-1}(|K_J|;\mathbb{Z})$$
\nd  with the convention $\widetilde{H}^{-1}(K_{\varnothing}; \mathbb{Z}) = \mathbb{Z}$.
\end{thm}
\nd  This leads to a formula for the bigraded Betti numbers in terms of the full subcomplexes
\begin{equation} \label{eqn:betti}
b^{-p,2q}(K) = \sum_{J\subset [m],  |J| =q}{\rm dim}\hspace{0.7mm}\widetilde{H}^{q-p-1}(|K_J|;\mathbb{Z}).
\end{equation}
>From Theorem \ref{thm:cohomofz}, (cf.~Theorem \ref{thm:baskakov}),  we conclude 
that the bigraded Betti numbers of the 
Tor module are the ranks of the contributions of $H^{\ast}(|K_J|)$ to the cohomology of $Z(K;(D^2,S^1))$. 
Further computations and properties of the bigraded and {\em multigraded\/} Betti numbers may be found in  \cite{aa2}
\cite{limon4} 
\cite{yu} and \cite{clu}. The calculations in \cite{clu} are used to verify the Halperin--Carlsson conjecture, (\cite[Remark 3]{clu}),
in the case of free torus actions on moment--angle complexes. This conjecture bounds the ranks of these tori by the
sum of the ranks of the cohomology groups.

An analogue of Theorem \ref{thm:cohomofz} is obtained  in \cite{lupan} for a generalization of moment--angle 
complexes arising from simplicial posets.  

The question whether the cohomology of a moment--angle manifold determines its diffeomorphism type is
currently an active area of research, \cite{bempp}.

\section{The exponentiation  property of polyhedral products}\label{sec:exp}
 A  product is defined on CW-pairs by
$$(X,A) \times (U,V) \;=\; \big(X\times U, (A\times U)\cup_{A\times V} (X\times V)\big).$$ 
Let $K$ be a simplicial complex on vertices $[m]$ and $\xa$ a family of CW-pairs. For a sequence of positive integers
$J = (j_1,j_2,\ldots,j_m)$, define a new family of pairs $\yb$ by
\begin{equation}\label{eqn:ybdefn}
(Y_i, B_i)  \;=\; (X_i,A_i)^{j_i}.
\end{equation}
It turns out that the combinatorics of a polyhedral product can detect this change. Let $K$ be a simplicial 
complex of dimension $n-1$ on vertices $\{v_1,v_2,\ldots,v_m\}$.
>From $K$  and a sequence $J$ as above, a new  simplicial complex $K(J)$ exists on $j_1 + j_2 + \cdots +j_m$ 
vertices, labelled
$$v_{11}, v_{12},\ldots,v_{1j_1},\;v_{21}, v_{22},\ldots,v_{2j_2},\;\;\ldots\;\; ,v_{m1}, v_{m2},\ldots,v_{mj_m}.$$
\nd It is characterized by the property that 
$$\big\{v_{i_{1}1},v_{i_{1}2},\ldots,v_{i_{1}j_{i_{1}}},
v_{i_{2}1},v_{i_{2}2},\ldots,v_{i_{2}j_{i_{2}}},\ldots,
v_{i_{k}1},v_{i_{k}2},\ldots,v_{i_{k}j_{i_{k}}}\big\}$$
\nd  is a  minimal non-face of $K(J)$ if and only if  $\{v_{i_1},v_{i_2},\ldots,v_{i_k}\}$ is a minimal 
non-face of $K$. Moreover, all minimal non-faces of $K(J)$ have this form.

An alternative explicit construction of the simplicial complex $K(J)$ may be found in either \cite[Construction $2.2$]{bbcg7}
or \cite{pb}. The construction of $K(J)$ is known variously in the literature as the {\em simplicial wedge construction\/} or
the {\em J-construction\/}.  \index{Wedge construction} \index{J-construction}

Next, we adopt the convention of denoting by $Z\big(K(J); (\underline{X}, \underline{A})\big)$ the polyhedral product
determined by the simplicial complex $K(J)$ and the family of pairs obtained from 
$(\underline{X}, \underline{A})$ by repeating each $(X_{i},A_{i})$, $j_{i}$ times in sequence.
The main structure theorem is the following.
\begin{thm}[\cite{bbcg7}]\label{thm:general.wedge}
Let $K$ be a simplicial complex with $m$ vertices and let $(\underline{X}, \underline{A})$ 
denote a family of CW-pairs and the family $\yb$ be defined as in \eqref{eqn:ybdefn}. Then 
$$Z\big(K(J);(\underline{X}, \underline{A})\big)\;=\;Z\big(K; (\underline{Y}, \underline{B})\big)$$
as subspaces of $X_1^{j_1} \times X_2^{j_2} \times \cdots \times X_m^{j_m}$.
\end{thm}

This theorem has applications  in toric topology and geometry. In particular, since we have
$(D^1, S^0)^2 \cong (D^2,S^1)$,  it implies that each moment-angle complex $Z(K;(D^2,S^1))$ is homeomorphic to 
a real moment-angle complex $Z(K(J);(D^1,S^0))$ where $J = (2,2,\ldots,2)$. 
Three other applications are discussed below.
\skp{0.2}
The cohomology of toric manifolds given by \eqref{eqn:zquotient} and \eqref{eqn:zequalszk} has a presentation 
determined by the combinatorics of $K$ and relations arising from the function $\lambda$. In cases when  $\mathcal{Z}$
in \eqref{eqn:zquotient} has the form $Z\big(K(J);(D^2,S^1)\big)$, for some $K$ and $J$, Theorem \ref{thm:general.wedge} 
can be used to give a simpler presentation of the cohomology of associated toric manifolds which uses the combinatorics 
of $K$ only. More details can be found in \cite{bbcg7}. 
\begin{rem}
A generalization of the Davis--Januszkiewicz construction \eqref{eqn:defn.tn} which accommodates
the construction here can be found in \cite{bbcg5}.
\end{rem}
An application to topological joins begins by defining a family of CW-pairs \index{Join, iterated}
$$\big(\underline{C(\divideontimes_{\!J}X}),\; \underline{\divideontimes_{\!J}X}\big)
\;\becomes\; \big\{\big(C(\underset{j_i}{\underbrace{X_i\ast X_i\ast\cdots\ast X_i}}),
\underset{j_i}{\underbrace{X_i\ast X_i\ast\cdots\ast X_i}}\big)\big\}_{i=1}^m.$$ 
for each sequence $J = (j_1,j_2,\ldots,j_m)$ and family of CW-complexes $\{X_i\}_{i=1}^{m}$.  Here, $C(-)$ denotes the 
unreduced cone.  It is shown in \cite{bbcg7} that the
equivalence $X\ast X \overset{\simeq}{\longrightarrow}  (CX \times X) \cup (X \times CX)$ iterates appropriately so that
the next theorem follows from Theorem \ref{thm:general.wedge}.
\begin{cor}[\cite{bbcg7}]
There is a homeomorphism of polyhedral products,
$$Z\Big(K; \big(\underline{C(\divideontimes_{\!J}X}),\; \underline{\divideontimes_{\!J}X}\big)\Big) 
\longrightarrow Z\big(K(J); (\underline{CX}, \underline{X})\big). $$              
\end{cor}
\begin{rem} Notice that in the case $X_i = S^1$ for all $i$, this corollary yields a homeomorphism
\begin{equation}\label{eqn:wedge}
Z\Big(K; \big(\underline{D}^{2J},\; \underline{S}^{2J-1}\big)\Big) 
\longrightarrow Z\big(K(J); (D^2, S^1)\big),
\end{equation} 
\nd where $\big(\underline{D}^{2J},\; \underline{S}^{2J-1}\big) = \big\{(D^{2j_i}, S^{2j_i - 1})\big\}_{i=1}^{m}$.
\end{rem}
\skp{0.2}
Theorem \ref{thm:general.wedge} also implies  an observation  about the action of the Steenrod algebra.
\begin{cor}[\cite{bbcg7, bbcg3}]\label{cor:steenrod}
There is an isomorphism of ungraded  rings 
$$H^{\ast}\big(Z(K;(D^2, S^1));\mathbb{Z}/2 \big) \longrightarrow 
H^{\ast}\big(Z\big(K(J); (D^2, S^1)\big);\mathbb{Z}/2 \big)$$
\nd which commutes with the action of the Steenrod algebra. 
\end{cor}
\begin{rem}
The construction of Theorem \ref{thm:general.wedge} has been generalized by Ayzenberg in \cite{aa2}; see also
\cite[Section 6]{bbcg5}.
\end{rem}
\section{Fibrations}\label{sec:early}
A direct extension of a result of Denham and Suciu \cite[Lemma $2.9$]{ds} yields the next theorem detailing the 
behaviour of polyhedral products with respect to fibrations. \index{Denham--Suciu fibration}
\begin{thm}\label{porter}
Let $p_i:(E_i,E'_i)\to(B_i,B_i')$ be a map of pairs, such that
$p:E_i\to B_i$ and $p|E'_i:E'_i\to B'_i$ are fibrations  over path connected CW complexes  with fibres
$F_i$ and $F'_i$ respectively. Let $K$ be a simplicial complex with
$m$ vertices.  If either $\underline{B} = \underline{B}'$ or $\underline{F} = \underline{F}'$, then the following
is a fibration:  
$$Z(K;(\underline{F},\underline{F}'))\to Z(K;(\underline{E},\underline{E}'))\to
Z(K;(\underline{B},\underline{B}').$$
\end{thm}
Special cases of this fibration were developed earlier by G. Porter, \cite{porter}. He proved that 
these fibrations exist in the special cases where $(X,A) = (X,\ast)$ with $A = \ast$ the base-point, and $K$ is the 
$q$-skeleton of the full $m$-simplex for $0 \leq q \leq m-1$.

 Natural consequences of Theorem \ref{porter} also arose earlier in the work
of T.~Ganea \cite{ganea} and A.~Kurosh \cite{kurosh}. G.~W.~Whitehead had defined a filtration
of the product $X_1\times\cdots \times X_m$, where the $j^{th}$-filtration is given by the space 
$$ W_j(X_1,X_2,\ldots,X_m)\;=\;\big\{(y_1,\ldots,y_m) \  |  \ y_i=\ast_i \in X_i, \ \hbox{for at least}\ m-j\
\hbox{values of}\ i\big\}. $$ 
\nd This space is the polyhedral product 
$$W_j(X_1,X_2,\ldots,X_m)\;=\;Z\big(\Delta[m-1]_{j-1};(\underline{X},\underline{\ast})\big),$$ 
\nd where $\Delta[m-1]_q$ denotes the $q$-skeleton of the  $(m-1)$-dimensional simplex. 
In this case, Porter's result gives an identification of the 
homotopy theoretic fibre of the inclusion 
$$W_j(X_1,X_2,\ldots,X_m) \;\subset\; \prod_{1 \leq i \leq m} X_i,$$  
which follows directly from the classical path-loop fibration over $X$ and
Theorem \ref{porter}. To see this, we shall extend to the right the fibration of pairs
$$\begin{CD}
(\Omega X,\Omega X)  @>{}>> (PX, \Omega X) @>{ e_1}>> (X,*),
\end{CD} $$
where $e_1$ is the map evaluating the endpoint of a path.
Consider  again  the cone $CY$ over a space $Y$. There is a map
$$\begin{array}{lccc}
\kappa \colon\hspace{-3mm}&X & \longrightarrow & X \times C(PX)\\[0.7mm]
                        &x           & \mapsto     &(x,[0,f_x])
\end{array}$$
\nd where $f_x\colon [0,1] \to X$ is the constant path satisfying $f_x(t) = x$. The map $\kappa$ is evidently
a homotopy equivalence.
Consider now $PX$ to be the subspace of $X \times C(PX)$ given by pairs $\big(f(1), [0,f]\big)$ for
$f \in PX$. This gives a pair $\big(X \times C(PX),PX\big)$.
Consequently, we get a fibration of pairs
$$\begin{CD}
\big(C(PX),\Omega X\big)  @>{}>> \big(X \times C(PX),PX\big) @>{\pi_X \times e_1}>>
(X,X)
\end{CD}$$
\nd for which $\pi_X:X \times C(PX) \to X$ is the natural projection. Applying 
Theorem \ref{porter}  yields the fibration
$$\begin{CD}
Z(K;(\underline{PX},\underline{\Omega X}))  @>{}>> Z(K;(\underline{X
\times PX},\underline{PX})) @>{\pi_X \times e_1}>>
Z(K;(\underline{X},\underline{X})).
\end{CD}$$
\nd Moreover, the equivalence of pairs $(X,\ast) \to \big(X\times C(PX),PX\big)$, gives the associated homotopy
fibration
\begin{equation}\label{eqn:porterfibration}
\begin{CD}
Z(K;(\underline{PX},\underline{\Omega X}))  @>{}>>
Z(K;(\underline{X},\underline{\ast})) @>{\pi_X \times e_1}>>
Z(K;(\underline{X},\underline{X})),
\end{CD}
\end{equation}
\nd which we record as the next corollary  which generalizes Porter's identification. 
\begin{cor} \cite[Corollary $2.32$]{bbcg2}  \label{cor:proof.of.porter}
If all of the $X_i$ are path-connected, the homotopy theoretic fibre
of the inclusion $ Z(K;(\underline{X},\underline{\ast})) \subset
Z(K;(\underline{X},\underline{X}))$ is
$Z(K;(\underline{PX},\underline{\Omega X})).$
In particular, the homotopy theoretic
fibre of the inclusion  $W_j(X_1,X_2,\ldots,X_m) \subset \prod\limits_{i=1}^{m}{X_i}$ 
is $$Z\big(\Delta[m-1]_{j-1};(\underline{PX},\underline{\Omega X})\big).$$
\end{cor}
\nd For a more general version of this result, see \cite[Proposition $5.1$]{prv}. See also
Section \ref{sec:aisapoint}, and for related information about generalized Whitehead products,
\cite[Section 7]{pt}, \cite{gt3} and \cite{ik6}, where the techniques used involve the {\em fat wedge filtration\/}.

The structure of the fibre $Z\big(\Delta[m-1]_{j-1};(\underline{PX},\underline{\Omega X})\big)$ is given by a second 
theorem of G.~Porter \cite{porter}.  After one suspension, this result also follows from the next theorem which gives
a decomposition of the suspension of 
$Z\big(\Delta[m-1]_{q};(\underline{CY},\underline{Y})\big)$, for a family of connected, pointed CW-complexes
$\{Y_i,\ast_i\big\}^m_{i=1}$.  Recall that
$\widehat{Y}^I$ is $Y_{i_1}\wedge\cdots\wedge Y_{i_k}$ if $I=(i_1,\ldots,i_k)$. 
\begin{thm}\label{T:1.18}
 Let $\{Y_i,\ast_i\big\}^m_{i=1}$ be a
family of connected, pointed CW-complexes.  Then there is a homotopy
equivalence
$$ \Sigma\big(Z\big(\Delta[m-1]_{q};(\underline{CY},\underline{Y})\big)\big)
\to\Sigma(\bigvee_{q+1<|I|\leq m}\bigvee_{t_I}(\Sigma
^{q+1}\widehat{Y}^I)),$$ 
where $t_I$ is the binomial coefficient
$\binom {|I|-1} {q+1}$.
\end{thm}
\nd   This result is deduced from \cite[Theorem $2.19$]{bbcg2}, by a counting argument 
which uses the fact that here, the full sub-complex $K_I$ is an $(|I|-1)$-simplex, and then enumerates
the spheres in the $q$-skeleton of the $(|I|-1)$-simplex. 

When $q=m-2$, $Z\big(\Delta[m-1]_{m-2};(\underline{X};\underline{\ast})\big)$
is the fat wedge  $W_{m-1}(X_1,X_2,\ldots,X_m)$  and the homotopy fibre is
$Z\big(\Delta[m-1]_{m-2};(\underline{PX},\underline{\Omega X})\big)$. In
this case, it follows that there is a homotopy equivalence
$$\Sigma\big(Z\big(\Delta[m-1]_{m-2};(\underline{PX},\underline{\Omega X})\big)\big) 
\to \Sigma(\Omega X_1 \ast\ldots \ast\Omega X_m).$$

T.~Ganea \cite{ganea}  identified the homotopy theoretic fibre of
the natural inclusion \\
$X \vee Y \to X \times Y$ as the join
$\Omega(X)\ast\Omega(Y)$ in the case $X$ and $Y$ are path-connected,
pointed spaces having the homotopy type of a CW-complex. This example is
a special case of the homotopy theoretic fibre for the inclusion of
a polyhedral product into a product corresponding to the
simplicial complex $K$ given by two  discrete  points.

\skp{1}
\section{Unstable and stable decompositions of the polyhedral product}\label{sec:splitting}
A Cartesian product $(X_1 \times \cdots \times X_m)$ of pointed, connected CW-complexes splits stably by 
a homotopy equivalence
$$H:\Sigma (X_1 \times \cdots \times X_m) \to \Sigma\big( \bigvee_{I \subseteq [m]} \widehat{X}^I\big)$$ where
$I$ runs over all the non-empty sub-sequences of $(1,2,\ldots,m)$. It is natural then to ask if a similar
splitting holds for the polyhedral product
$$Z(K;\xa)\; \subset \; X_1 \times \cdots \times X_m.$$
\nd In other words, does a stable splitting of the polyhedral product exist for simplicial complexes $K \neq \Delta^{m}$\/?
Indeed, in the setting of moment--angle complexes, such a splitting was suggested by the cohomology result of 
V.~Baskakov \cite{bask} following. \index{Baskakov's Theorem}
\begin{thm}\label{thm:baskakov}
There is an isomorphism of rings
$$H^k\big(Z_K;(D^2,S^1)\big) \; \cong\; \bigoplus_{I \subset [m]}\widetilde{H}^{k-|I|-1}(|K_I|),$$
\nd where the product on the right hand side is that described in Section \ref{sec:cohomology}.
\end{thm}
\vspace{-0.1\baselineskip}
\nd   This result is deduced from those in Section \ref{sec:cohomofamac}.
 This decomposition is related to a splitting of the cohomology of the complements
of subspace arrangements as follows. \index{Subspace arrangements}  A simplicial complex $K$ determines a complex
coordinate arrangement of subspaces whose  complement is homeomorphic to $Z(K;(\mathbb{C},\mathbb{C^\ast}))$. 
The moment--angle complex $Z(K;(D^2,S^1))$ is a strong deformation retract of $Z(K;(\mathbb{C},\mathbb{C^\ast}))$,  
by a result of N.~Strickland  \cite[Proposition 20]{strickland}, (cf.~V.~Buchstaber and T.~Panov \cite[Theorem $5.2.5$]{bp3}).
In \cite{gm}, M.~Goresky and  R.~MacPherson derive a decomposition for the cohomology of the complements
of subspace arrangements which gives a cohomology splitting analogous to Theorem \ref{thm:baskakov}. 
A more direct proof was given by G.~Ziegler and R.~\v{Z}ivaljevi\'{c} in \cite{zz}. 
\skp{0.1}
The first geometric result in this direction is the {\em unstable\/} splitting due to J.~Grbi\'c and S.~Theriault, \cite[Theorem 1]{gt1}. 
\index{Moment--angle complex, unstable splitting}
They show that when
$K$  consists of $m$ discrete  points, then there is a homotopy equivalence
\begin{equation}\label{eqn:disjointpoints}
Z(K;(D^2,S^1)) \simeq Z(K;(\mathbb{C},\mathbb{C^\ast})) \longrightarrow \bigvee_{k=2}^{m}(k-1)\binom{m}{k}S^{k+1}.
\end{equation}
\nd Their proof uses an inductive argument based on a version of Corollary \ref{cor:proof.of.porter}. As noted above, the space
$Z(K;(\mathbb{C},\mathbb{C^\ast}))$ corresponds to the complement of a codimension-two coordinate subspace arrangement.
Later in \cite{gtshifted}, they extended this result to a class of simplicial complexes called {\em shifted}, 
\index{Simplicial complex, shifted} and to simplicial complexes
which can be obtained from shifted complexes by certain elementary topological operations.
\begin{defin}\label{def:shifted}  \index{Simplicial complex, shifted}
A simplicial complex $K$ is called {\em shifted\/} if there is an ordering of its vertices such that  whenever $v$ is a vertex
of $\sigma \in K$  and $v' < v$, $(\sigma \smallsetminus v) \cup v' \in K$.
\end{defin}
\begin{thm}\cite[Theorem 1.1]{gtshifted}\label{thm:gtshifted} If $K$ is a shifted complex, then $Z(K;(D^2,S^1))$ is homotopy 
equivalent to a wedge of spheres.
\end{thm}
Following a conjecture made by the authors in \cite{bbcg2}, this result has been extended in \cite{gt2} and \cite{ikshifted}
to the case $Z(K;(\underline{CX},\underline{X}))$, as described in Theorem  \ref{thm:cxcommax}  below. Indeed, under a 
condition on $K$, known as  {\em dual sequentially 
Cohen-Macaulay over} $\mathbb{Z}$, there is the following improvement  over previous results, including those of \cite{gw}.
 (Recognizing that shifted implies dual shifted, arguments from \cite{bwachs} and \cite{bww} can be amalgamated to deduce 
that $K$ shifted implies that it is dual sequentially Cohen-Macaulay over $\mathbb{Z}$.) 
\begin{thm}\cite[Theorems 1.3 and 1.4]{ik5}\label{thm:cxcommax}  \index{polyhedral product, unstable splitting}
If $K$ is dual sequentially Cohen-Macaulay over $\mathbb{Z}$, then
\begin{equation}\label{eqn:cxx}
Z(K;(\underline{CX},\underline{X})) \simeq \bigvee_{\varnothing \neq I \subset [m]}\Sigma{|K_{I}|}\wedge \widehat{X}^I.
\end{equation}
\end{thm}
\nd  Notice that if $X = S^1$, the right hand side is a wedge of spheres. 
\begin{exm}
A simple example of the splitting in Theorem \ref{thm:cxcommax} occurs for 
the simplicial complex consisting of two  discrete  points, $K = \big\{\varnothing , \{v_1\}, \{v_2\}\big\}$. Here $m = 2$ and 
$$(\underline{CX},\underline{X}) = \{(CX_1,X_1), (CX_2,X_2)\}.$$
In this case, \eqref{eqn:cxx} sees the well known homotopy equivalence 
$\;X_{1}\ast X_{2}\; \simeq\;  \Sigma{(X_{1}\wedge X_{2})}\;$
in the form
$$CX_{1}\times X_{2} \cup_{X_{1}\times X_{2}} X_{1}\times CX_{2}\; \simeq \;
CX_{1} \vee CX_{2}\vee \Sigma{(X_{1}\wedge X_{2})} .$$
\end{exm}
For general $\xa$ there is the stable splitting
of \cite{bbcg1} and \cite{bbcg2}. \index{Polyhedral product, stable splitting}
\begin{thm}\cite[Theorem 2.10]{bbcg2}\label{thm:bbcgsplitting}
Let $K$ be an abstract simplicial complex  on vertices $[m]$.  Given $\xa) = \{(X_i , A_i)\}_{i=1}^m$ where $(X_i , A_i , x_i )$ 
are pointed  pairs  of CW-complexes.
then there is a natural pointed homotopy equivalence 
\begin{equation}\label{eqn:splitting}
H\colon \Sigma{\big(Z\big(K;\xa\big)\big)} \longrightarrow 
\Sigma\Big(\bigvee_{I\subseteq [m]}\widehat{Z}\big(K_I; (\underline{X}, \underline{A})_I\big)\Big)
\end{equation}
where $(\underline{X}, \underline{A})_I$ denotes the restricted family of CW-pairs $\big\{(X_i,A_i\big)\}_{i\in I}$.
\end{thm}
\begin{rem}
The spaces $Z\big(K;\xa\big)$ generally do not decompose as a wedge before suspending. One example is
the simplicial complex $K$ determined by a square, with $4$ vertices and $4$ edges, for which 
$Z\big(K;(D^2,S^1)\big)$ is $S^3 \times S^3$. 
\end{rem}
The homotopy equivalence \eqref{eqn:splitting} is induced from the following decomposition which is well known, 
\cite{milnor, james}. 
\begin{thm} \label{thm:milnor}
Let $(Y_i,y_i)$ be pointed CW-complexes. There is a pointed, natural
homotopy equivalence $$H:\Sigma (Y_1 \times \cdots \times Y_m) \to
\Sigma\big( \bigvee_{I \subseteq [m]} \widehat{Y}^I \big)$$ where $I$ runs over
all the non--empty subsequences of $(1,2,\ldots,m)$. Furthermore,
the map $H$ commutes with  colimits.
\end{thm}
\nd  It is shown in \cite{bbcg2} that the map $H$ of Theorem \ref{thm:milnor} induces an equivalence between 
the diagrams which define the colimits on either side of \eqref{eqn:splitting}, yielding Theorem \ref{thm:bbcgsplitting}.
\skp{0.2}
In many important cases, the spaces 
which appear on the right hand side of  \eqref{eqn:splitting} can be identified explicitly. Below are a few of the most common
of these. We begin with a definition.
\begin{defin}\label{defn:link}
For $\sigma$  a simplex in $K$, $\text{lk}_{\sigma}(K)$  {\em the link of\/} $\sigma$ {\em in\/} $K$,
is defined to be the simplicial complex  for which
$$ \tau \in \text{lk}_{\sigma}(K)\quad \text{if and only if}\quad \tau \cup \sigma \in K.$$
\end{defin}
\begin{thm}\cite{bbcg2}\label{thm:nullincl}
Let $K$ be an abstract simplicial complex  on vertices $[m]$  and $I \subseteq [m]$. 
If the family $(\underline{X},\underline{A})$
has the additional property  that the inclusion $A_i\subset X_i$ is null-homotopic for all $i$, then there
is a homotopy equivalence
\begin{equation}\label{eqn:ainull}
\widehat{Z}\big(K_I;(\underline{X}, \underline{A})_I\big)\longrightarrow \bigvee\limits_{\sigma\in
K_I} |{\rm{lk}}_{\sigma}(K_I)|\ast\widehat{D}(\sigma)
\end{equation}
\nd which, when combined with \eqref{eqn:splitting}, gives a homotopy equivalence
$$\Sigma{\big(Z\big(K;\xa\big)\big)} \longrightarrow 
\Sigma\Big(\bigvee\limits_{I\subseteq [m]}\big(\bigvee\limits_{\sigma\in
K_I} |{\rm{lk}}_{\sigma}(K_I)|\ast\widehat{D}(\sigma)\big)\Big).$$
\end{thm}
\nd  (Here, $\widehat{D}(\sigma)$ is defined in Section \ref{sec:Introduction}.) 
The homotopy equivalence \eqref{eqn:ainull} is a generalization of the wedge lemma of 
Welker-Ziegler-\v{Z}ivaljevi\'c in \cite{wzz}, applied to the smash polyhedral product. 
 
Next, consider all homology to be taken with  coefficients in a field  $k$. Define the Poincar\'e series 
of a  pointed  space $X$ of finite type  by 
$$P(X,t) = \sum d_n(X)t^n$$ 
where $d_n(X) = \mbox{dimension}_{\mathbb F}H_n(X; \mathbb F),$
and the {\em reduced\/} Poincar\'e series by 
$$\overline{P}(X,t) = -1 + P(X,t).$$ This series behaves well with respect to 
wedges and smash products. Hence, (for example),  for pairs satisfying the hypotheses of Theorem \ref{thm:nullincl} we have
the next result.
\index{Polyhedral product, Poincar\'e series}

\begin{lem} \label{lem:Poincare.series.null}
Assume that homology is taken with field coefficients $\mathbb F$ and that path connected pairs of pointed CW-complexes 
$(\underline{X},\underline{A})$,  of finite type, have the property that the
inclusion $A_i\subset X_i$ is null-homotopic for all $i$. Then  
$$\overline{P}\big(\widehat{Z}(K;(\underline{X},\underline{A})),t\big) =   \sum_{ \sigma\in K}    
t\overline{P}\big(|{\rm{lk}}_{\sigma}(K_I)|,t\big)\overline{P}\big(\widehat{D}(\sigma) ,t\big).$$ 
\end{lem}  
 Notice that if in Theorem \ref{thm:nullincl}  all  the spaces $A_i$ are contractible, then the spaces
$\widehat{D}(\sigma)$  are contractible unless $I$ represents a simplex $\sigma$ in $K$, that is $K_I$ is an
$(|I|-1)$--simplex; in this case, $\widehat{D}(\sigma) = \widehat{X}^I$ and the next theorem follows.
\begin{thm}[\cite{bbcg2}] \label{thm:more.contractible.A} 
Let $K$ be an abstract simplicial complex with $m$ vertices and
$(\underline{X},\underline{A})$ has the property that all the $A_i$
are contractible. Then there is a homotopy equivalence
$$\widehat{Z}\big(K_I;(\underline{X}, \underline{A})_I\big)\longrightarrow  \widehat{X}^{I}$$
yielding a homotopy equivalence
\begin{equation}\label{eqn:more.contractible.A}
\Sigma\bigl(Z\bigl(K;(\underline{X},\underline{A})\bigr)\bigr) \longrightarrow  \Sigma( \bigvee_{I \in K} \widehat{X}^{I}).
\end{equation}
\end{thm}
The next example is the one for which all the spaces $X_i$ are contractible. Now, for each $I \subseteq [m]$, the spaces
$\widehat{D}(\sigma)$  in \eqref{eqn:ainull} are contractible with the possible exception of the case $\sigma = \varnothing$,
the empty simplex, which gives $\widehat{D}(\sigma) = \widehat{A}^I$. In this case $|{\rm{lk}}_{\varnothing}(K_I)| \simeq |K_I|$
and \eqref{eqn:ainull} simplifies to yield the next theorem.
\begin{thm}\label{thm:xicontractible} \cite{bbcg2}
Let $K$ be an abstract simplicial complex with $m$ vertices and suppose that the collection
$(\underline{X},\underline{A})$ has the property that all the $X_i$
are contractible. Then there  are homotopy equivalences
$$\widehat{Z}\big(K_I;(\underline{X}, \underline{A})_I \big)\longrightarrow |K_I|\ast\widehat{A}^{I}
\longrightarrow \Sigma (|K_I|\wedge\widehat{A}^{I})$$
yielding a homotopy equivalence
$$\Sigma\big(Z(K; (\underline{X},\underline{A}))\big) \longrightarrow
\Sigma \big(\bigvee_{I\notin K} |K_I|*\widehat{A}^I\big).$$
\end{thm}
\nd The last statement of this theorem follows on observing that if $I \in K$ then $|K_I|$ is contractible. 
 The case $(X_i, A_i, x_i) =  (D^{n+1},S^n, \ast)$ for all $i$, $n \geq 1$, yields the next corollary. 
\begin{cor}\label{cor:bigmac}
There are homotopy equivalences
$$\Sigma \big(Z(K;(D^{n+1},S^n))\big)\to
\Sigma \big(\bigvee_{I \notin K} |K_I|*S^{n|I|}\big) \to \bigvee_{I \notin
K} \Sigma^{2 + n|I|}|K_I|.$$
\end{cor}
In the case $n=1$, which corresponds to the case of complements of complex coordinate subspace arrangements,
this corollary implies the Goresky--MacPherson cohomological splitting mentioned earlier. To observe this, the reduced 
cohomology of full sub-complexes of $K$ needs to be related to the reduced homology of certain links in the Alexander 
dual of $K$. This is done in, for example, \cite[Corollary $2.4.6$]{bp2}. 
\section{Equivariance of the stable splitting and an application to number theory}\label{sec:equivariance}
 In the case that all the pairs $(X_i,A_i)$ are the same,  the automorphism group of $K$, 
{\em Aut}$(K) \subseteq \Sigma_m$, acts in a natural way on both sides of 
the splitting \eqref{eqn:splitting}.   Using James--Hopf invariants, A.~Al-Raisi \cite{ar} shows that there is a stable 
splitting, as in \eqref{eqn:splitting}, which is equivariant with respect to the action of {\em Aut}$(K)$.  
\index{Stable splitting, equivariance}
\begin{thm}[\cite{ar}]\label{thm:alraisi}
\nd  The adjoint of the stable 
decomposition \eqref{eqn:splitting}, 
$$\theta_H: Z(K;(X,A))\; \longrightarrow\; J(\bigvee_{I \subseteq [m]}\widehat{Z}(K_I;(X,A)_I)) \longrightarrow\;
\Omega\Sigma\Big(\bigvee_{I\subseteq [m]}\widehat{Z}\big(K_I; (X,A)_I\big)\Big)\; $$ 
\nd is $Aut(K)$-equivariant. 
\end{thm}
 
 When $K = K_n$ is the boundary dual of an $n$-gon,  (see Section \ref{sec;bpapproach}),  a theorem of Coxeter 
 identifies the real moment--angle manifold  $Z(K_n ; (D^1,S^0))$.
\begin{thm}[\cite{coxeter}]
The real moment--angle manifold $Z(K_n ; (D^1,S^0))$  corresponding to an  $n$-gon, $n \geq 4$, is an 
oriented surface of genus 
$g=1+(n-4)2^{n-3}$.
\end{thm}
\begin{rem} 
The simplicial complex $K_n$ has also been used  by S.~Das \cite{das}, to show directly that the genus of the  
{\em hypercube graph\/}, the one-skeleton of the $n$-cube, is $1+(n-4)2^{n-3}$. This is a result due  originally to
G.~Ringel \cite{ringel} and L.~Beineke and F.~Harary \cite{bhar}. The computation is done by embedding equivariantly the 
hypercube graph into the real moment-angle complex $Z(K_{n}; (D^1,S^0))$. 
\end{rem}

 In this simple case, where {\em Aut}$(K_n) = D_{2n}$, the dihedral group with cyclic subgroup $C_n$, Theorem
\ref{thm:alraisi} has interesting consequences. 
The real moment--angle manifold is the conduit through which A.~Al-Raisi uses Theorem \ref{thm:alraisi} to link combinatorics 
of $K_n$
and the action of $C_n$ on a choice of basis for the cohomology of  a surface of genus  $g=1+(n-4)2^{n-3}$.  

Briefly, one first counts the number of orbits of the action of $C_n$ by considering, up to rotation, the boundary of an 
$n$-gon as a $2$-colored beaded necklace  with  $n$ beads. Certain colorings are shown to 
correspond to generators of  $H^{\ast} (Z(K_n ; (D^1,S^0)))$.  A theorem of Burnside gives the number of orbits of the action 
of $C_n$ arising from this choice of basis as
$$ \dfrac{1}{n} \Big( \underset{d|n}{\sum} \phi(d) 2^{n/d} -(n+1)\Big)$$ 
where $\phi(d)$ is the Euler phi function.

Next,  we think of the sequence of beads as a word made from two letters. A word of length $n$ is called an {\em aperiodic word\/}
if it has $n$ distinct {\em rotations\/}.  Here a rotation of a word $W = a_1a_2\ldots a_n$ is given by
$r(W) = a_na_1a_2\ldots a_{n-1}$. An equivalence class of an aperiodic word under rotation is called a {\em primitive necklace.\/}
A.~Al-Raisi shows that the rank of $H_{1} (Z(K_n ; (D^1,S^0)))$ is given in terms of the number of primitive necklaces.
The number of primitive $n$-necklaces on an alphabet of size $k$, denoted by $M(k,n)$, was computed by Moreau in the 
nineteenth century, \cite{moreau}
$$M(k,n) = \dfrac{1}{n}\underset{d|n}{\sum}\mu(d)k^{\frac{n}{d}}$$
where $\mu$ is the M\"{o}bius inversion function. 
The number of orbits of size $d$,  where $d|n$ and $ 1<d \leq n$, can be counted with the aid of this theorem.
The number of orbits with  exactly $n$ elements is given by 
$$ \dfrac{1}{n} \Big( \underset{d|n}{\sum} \mu(d) 2^{n/d}-(n-1) \Big).$$
The number of orbits of size $d<n$, where $d|n$, is given by
$$\dfrac{1}{d} \underset{d_1|d}{\sum} \mu(d_1)2^{d/d_1}.$$
Finally, equating the two different orbit counts,  we arrive at Al-Raisi's proof  of the cyclotomic identity below. 
 \index{Cyclotomic identity}
\begin{thm}[\cite{ar}]
For $n \geq 4$ there is an identity
$$\dfrac{1}{n} \Big( \underset{d|n}{\sum} \mu(d) 2^{n/d}\Big)-(n-1)  +  
\underset{1<d<n}{\underset{d|n}{\sum}} \Big(\dfrac{1}{d}  \underset{d_1|d}{\sum} \mu(d_1)2^{d/d_1}\Big) \;=\; 
\dfrac{1}{n} \Big( \underset{d|n}{\sum} \phi(d) 2^{n/d} \Big) -(n+1) $$
\nd relating the M\"obius function $\mu(d)$ and the Euler phi function $\phi(d)$.
\end{thm}

The rank of $H_1$ above is incidentally given by the rank of certain Lie tensors
in a free Lie algebra as first computed by E.~Witt  \cite{witt}. Moreau's formula and Witt's formula agree with one another. 
This connection between polyhedral products and free Lie algebras is as yet not clearly understood. 

 These structures together with connections to Hochschild homology as well as new filtrations of classifying spaces 
will appear in \cite{acv}. 
Additional interesting work subsequent to \cite{ar}, but from a different point of view, was developed by X.~Fu and 
J.~Grbi\'c.  They construct a sequence of simplicial complexes $\{L_m\}$ with each $L_m$ obtained from the 
$m$-cube by a vertex cut, \cite[Construction $4.4$]{fg} and verify a form of $\Sigma_m$--representation stabilty
in characteristic zero, for the homology of $Z(L_{m}; (X,A))$. 

In a recent preprint, S.~Cho, S.~Choi and S.~Kaji \cite{cck}, also 
address the representations afforded by this action. They compute some 
actions and the associated representations on the homology of $Z(K_n ; (D^1,S^0))$ in characteristic zero.

\section{The case that $A_i = \ast$ for all $i$}\label{sec:aisapoint}
The stable structure of the polyhedral product \apoint \;is given by Theorem \ref{thm:more.contractible.A}.
The importance of this particular case of the polyhedral product first came to light in the development of
toric topology with the following observation.
\begin{thm}\cite[Theorem $4.3.2$]{bp2}.\label{thm:djequalsbp}
The  inclusion $i\colon Z\big(K;(\mathbb{C}{\rm P}^\infty, \ast)\big) \longrightarrow (\mathbb{C}{\rm P}^\infty)^m$\\
factors into a composition of a homotopy equivalence
\begin{equation}\label{eqn:djequalsbp}
Z\big(K;(\mathbb{C}{\rm P}^\infty, \ast)\big) \; \longrightarrow \; ET^m \times_{T^m} Z(K;(D^2,S^1))
\end{equation}
and the fibration $ET^m \times_{T^m} Z(K;(D^2,S^1)) \longrightarrow BT^m = (\mathbb{C}{\rm P}^\infty)^m.$  
The action of the torus $T^m$ here is that of \eqref{eqn:cubical}. 
\end{thm}
This gives one of the most important  fibrations  in this subject;
\begin{equation}\label{eqn:basicfib}
Z(K;(D^2,S^1)) \stackrel{\widetilde{\omega}}{\longrightarrow} Z(K;(\mathbb{C}{\rm P}^\infty,\ast)) \longrightarrow 
(\mathbb{C}{\rm P}^\infty)^m.
\end{equation}
\nd The cohomology of the right-hand side of \eqref{eqn:djequalsbp} was determined by geometers in the context
of toric geometry, \cite{brion, bcp}. In the setting of toric topology, it appears in the paper of 
M.~Davis and T.~Januszkiewicz, \cite[Theorem $4.8$]{dj}. This calculation, together with the theorem, implies the following 
result in algebraic combinatorics, for which independent proofs
exist.
\begin{cor}[\cite{bp3,bbcg2}] \label{cor:sr} 
Let $K$ be a simplicial complex. Then
$$H^\ast\big(Z\big(K;(\mathbb{C}{\rm P}^\infty, \ast)\big);\mathbb{Z}\big) \; \cong\;  \mathbb{Z}(K)$$
\nd the Stanley--Reisner ring of $K$.
\end{cor}
 
Recall that Theorem \ref{thm:cohomofz} tells us that $\mathbb{Z}(K)$, determines the cohomology of 
the moment--angle complex $Z(K;(D^2,S^1))$.

An analogue of Corollary \ref{cor:sr}  holds over any ring $R$ for the algebra $R(K)$. Somewhat surprisingly, 
$R(K)$ determines $K$ by a 
theorem of  W.~Bruns and  J.~Gubeladze,  \cite{bg}, \cite[Theorem 5.27 and Example $5.28$]{bgbook}.   
\begin{thm}[\cite{bg}]\label{thm:commcase} Let $K_1$ and $K_2$ be simplicial complexes satisfying 
$R(K_1) \cong R(K_2)$ as $R$-algebras.
Then $K_1$ and $K_2$ are isomorphic as simplicial complexes. 
\end{thm}  
Over any field $k$ there is also  the {\em exterior\/} Stanley--Reisner ring \index{Stanley--Reisner ring, exterior} 
over $k$, which we denote by $\Lambda_{k}(K)$. 
It is the quotient of a graded exterior algebra on $m$ one-dimensional generators, by the same two-sided ideal as in 
\eqref{eqn:srring}. In a manner entirely analogous to the proof of Corollary \ref{cor:sr} from 
Theorem \ref{thm:more.contractible.A}, 
we have 
\begin{equation}\label{eqn:stretch}
\Lambda_{k}(K) \cong H^{\ast}\big(Z(K; (S^1, \ast)) ; k\big).
\end{equation}
 A proof that \eqref{eqn:stretch} also determines $K$ up to isomorphism, follows by an argument similar to that for
Theorem \ref{thm:commcase}, see \cite[Exercise $5.12$]{bgbook}. Independently, 
using the theory of algebraic groups,  C.~Stretch, (unpublished), asserts the same result. 
In order to compare the Stanley--Reisner ring to {\eqref{eqn:stretch} from a geometric perspective, it is convenient to work over
$k = \mathbb{Z}_{2}$ and to use the ungraded isomorphism of rings, (Theorem \ref{thm:more.contractible.A}, again),
\index{Loop spaces} 
\begin{equation}\label{eqn:ungradedextsr}
\Lambda_{\mathbb{Z}_{2}}(K)\cong H^{\ast}\big( Z(K; (S^2,\ast));\mathbb{Z}_{2}\big), \quad \text{(as ungraded rings)}.
\end{equation} 
The   fibrations associated to \eqref{eqn:porterfibration} or \eqref{eqn:coneloopfib}, give a 
commutative diagram
\begin{equation}\label{eqn:srandext}
\begin{CD}
Z(K;(D^2,S^1))  @>{}>>  Z(K;(\mathbb{C}{\rm P}^\infty, \ast))    @>{}>>  (\mathbb{C}{\rm P}^\infty)^m \\
@AA{\delta}A       @A{i}AA  @AA{i}A \\
Z(K;(C(\Omega{S^2}), \Omega{S^2}))@>{}>> Z(K;(S^2,\ast) @>{}>> (S^2)^m
\end{CD}
\end{equation}
where the right hand vertical maps are determined by the inclusion $S^2 \xrightarrow{i} \mathbb{C}{\rm P}^\infty$. In cohomology
with $\mathbb{Z}_{2}$  coefficients,
the middle map  determines the natural map from  the  ring $\mathbb{Z}_{2}(K)$ to the algebra
$\Lambda_{\mathbb{Z}_{2}}(K)$ of \eqref{eqn:ungradedextsr}. The map 
$$\Omega{S^2} \stackrel{\Omega{i}}{\longrightarrow} \Omega{\mathbb{C}{\rm P}^\infty}\;\simeq\; \Omega{BS^1} \;\simeq\; S^1$$
induces the vertical map $\delta$, and is given by the projection
$\Omega{S^2} \;\simeq\; S^1\times \Omega{S^3} \longrightarrow S^1$,
since this  detects  the generator of $\pi_{1}(\Omega{S^2})$. 
\skp{0.4}
\nd\centerline{\LARGE{$\longrightarrow\!\longleftarrow$}}
\skp{0.1}
 It follows from \eqref{eqn:djequalsbp} that the Borel construction is a CW-subcomplex of $ (\mathbb{C}{\rm P}^\infty)^m$. 
In particular it has cells in even degree only. 
In the literature the space $Z\big(K;(\mathbb{C}{\rm P}^\infty, \ast)\big)$  is often denoted by $DJ(K)$ or $DJ_K$ and is
called the {\em Davis--Januszkiewicz space\/}. \index{Davis--Januszkiewicz space}

The integral and rational formality of the Davis--Januszkiewicz spaces is studied  by  M.~Franz in \cite{franz2}, 
and further by  D.~Notbohm and N.~Ray in \cite{nr}, using model category theory. They show that the rationalization of $DJ(K)$ is 
unique for a class of simplicial complexes 
$K$. 

\begin{rem}  Under suitable freeness conditions, there is an analogue, 
$H^\ast\big(Z\big(K;(X, \ast)\big);\mathbb{Z}\big)$, 
of the Stanley--Reisner ring for any space $X$ and simplicial complex $K$, which 
specializes to the classsical Stanley--Reisner ring in case $X$ is $\mathbb{C}{\rm P}^\infty$. The result follows from 
Theorem \ref{thm:more.contractible.A}, and for more details, \cite[Theorem $2.35$]{bbcg1}. Example \ref{exm:sr} in 
Section \ref{sec:cohomology} gives a different perspective and a more concise description of the ring under discussion 
here. It is in fact the case that {\em all\/} monomial ideal rings can be realized as  cohomology rings of spaces  related to 
 certain polyhedral products, \cite{monomials}. 

Unstable phenomena in the homotopy of the 
left-hand side of \eqref{eqn:djequalsbp} has been investigated by D.~Allen and J.~La Luz \cite{la1,al1} using the Unstable 
Novikov Spectral Sequence. In particular, their study of the derived functors of the indecomposables of the 
Stanley--Reisner ring, leads to an interesting and mysterious invariant of simplicial complexes. 
\end{rem}

Next, we consider the relationship between the space 
$Z(K;(S^1,\ast))$ and the right--angled Artin group determined by $K$.
\begin{defin}\label{defn:graph.product}  \index{Graph product}
Given a simplicial graph $\Gamma$ with vertex set $S$ and a family of groups $\{G_s\}_{s\in S}$, their graph product
$\prod_{\Gamma}G_s$ is the quotient of the free product of the groups $G_s$ by the relations that elements of $G_s$ 
and $G_t$ commute whenever $\{s,t\}$ is an edge of $\Gamma$. In the case that $G_i  = \mathbb{Z}$ and 
$\Gamma = SK$, the one-skeleton
of $K$, the graph product is called the {\em right--angled Artin group\/}, $RA_K$, corresponding to $K$.   When 
$G_i = \mathbb{Z}_2$, the graph product is called the {\em right--angled Coxeter group\/}, $RC_K$. 
\index{Right--angled Artin group} \index{Right--angled Coexeter group} A simplicial complex for which every collection of 
pairwise adjacent vertices spans a 
simplex, is called a {\em flag complex\/}.  \index{Flag complex}
\end{defin}
\begin{thm}[\cite{kb,kr,do,chardav,pv}]\label{thm:kpi1} \index{Polyhedral product, aspherical}
The space $Z(K;(S^1,\ast))$ is aspherical if and only if $K$ is a flag complex, in which case, 
$\pi_1\big(Z(K;(S^1,\ast))\big) = RA_K$.
\end{thm}
Related results, \cite{davis,pv}, identify $\pi_{1}\big(Z(K;(\mathbb{R}{\rm P}^\infty,\ast))\big)$ as a right--angled 
Coxeter groups $RC_K$. From this follows the fact that $\pi_{1}\big(Z(K;(D^1,S^0))\big ) = [RC_K,RC_K]$, the commutator
subgroup.

A general result in this direction is due to M.~Stafa.  A proof is sketched below using
a recent result about the  homotopy type of $\Omega{Z}(K;(\underline{X},\underline{\ast}))$ due to T.~Panov and 
S.~Theriault, \cite{pt}.
\begin{thm}[\cite{ms2}]\label{thm:general.kpi1}
Let $G_1,G_2,\ldots,G_m$ be non--trivial discrete groups and $K$ be a simplicial complex with $m$ vertices. Then 
$Z(K;(\underline{BG},\underline{\ast}))$ is an Eilenberg--Mac Lane space if and only if $K$ is a flag complex. Equivalently, 
$Z(K;(\underline{EG},\underline{G}))$ is 
an Eilenberg--Mac Lane space if and only if $K$ is a flag complex.
\end{thm}
\nd {\em Proof\/}.~({\sc sketch.})
For the simplicial complex $V_K$, consisting of the vertices of $K$, there is a natural map
$$X_1\vee X_2\vee\cdots\vee X_m \;\xrightarrow {\simeq} \;Z(V_K; (\underline{X},  \underline{\ast}))\;
\xrightarrow{i} \;Z(K; (\underline{X},  \underline{\ast}))$$ 
where the first equivalence is that of \eqref{eqn:more.contractible.A} and the second map is induced by the inclusion 
$V_K \hookrightarrow K$.  For a flag complex $K$, T.~Panov and S.~Theriault show in \cite{pt} that $\Omega{i}$ has a right 
homotopy inverse
\begin{equation}\label{eqn:pt}
\Omega(X_1\vee X_2\vee\cdots\vee X_m)\; {\longleftarrow}\;\Omega Z\big(K; (\underline{X},  \underline{\ast})\big).
 \end{equation}
 \nd so that $\Omega Z\big(K; (\underline{X},  \underline{\ast})\big)$ is a retract of $\Omega(X_1\vee X_2\vee\cdots\vee X_m)$. 
\skp{0.1}
When $X_i = BG_i$ it is  known that $BG_1 \vee BG_2 \vee\cdots\vee BG_m \simeq B(G_1\ast G_2 \ast\cdots\ast G_m)$,
and so for $K$ a flag complex, it follows immediately that 
$Z(K;(\underline{BG},\underline{\ast}))$ is an Eilenberg--Mac Lane space.

If $K$ is not a flag complex, then $\partial\Delta^n$, the boundary of an $n$--simplex,  is a full subcomplex of  
$K$ for some $n>1$.  In this case,  \cite[Proposition $3.3.1$]{ds} implies that
$Z\big(\partial\Delta^n; (\underline{BG},\underline{\ast})\big)$ splits off   from
$Z\big(K; (\underline{BG},\underline{\ast})\big)$. Hence it suffices to show that 
$Z\big(\partial\Delta^n; (\underline{BG},\underline{\ast})\big)$ has non-trivial higher homotopy groups.
\nd For the simplicial complex $\partial\Delta^n$, the fibration of Theorem \ref{porter} takes the form
\begin{equation}\label{eqn:dsforbg}
Z\big(\partial\Delta^n; (\underline{EG},\underline{G})\big) \longrightarrow 
Z\big(\partial\Delta^n; (\underline{BG}), \underline{\ast})\big) \longrightarrow
BG_1 \times BG_2 \times \cdots\times BG_{n+1}.
\end{equation}
Next, since $\partial\Delta^n$ is a shifted simplicial complex, (Definition \ref{def:shifted}), and $EG$ is contractible,
Theorem \ref{thm:cxcommax} gives
$$ Z\big(\partial\Delta^n; (\underline{EG},\underline{G})\big) \simeq 
S^n \wedge G_1 \wedge G_2 \wedge \cdots \wedge G_{n+1}.$$
\nd Since the groups $G_i$ are discrete for all $i$, $ Z\big(\partial\Delta^n; (\underline{EG},\underline{G})\big)$ 
becomes a wedge of $n$-spheres. The assumption $n>1$, the Hilton-Milnor theorem, and the homotopy sequence of
\eqref{eqn:dsforbg} imply that
$Z\big(\partial\Delta^n;(\underline{BG},\underline{\ast})\big)$ has non-trivial higher homotopy groups which completes the proof.

M.~Stafa \cite{ms} also uses polyhedral  products to construct  monodromy representations of a finite product of discrete
groups into outer automorphisms of free groups.  We close this section by noting that  the classifying spaces of certain 
groups derived from right-angled Coxeter groups through quandles, are described in terms of polyhedral products, with 
$X_i \neq \ast$ and $A_i \neq \ast$, by D.~Kishimoto in \cite{kishimoto}.

\section{The cohomology of polyhedral products and a spectral sequence}\label{sec:cohomology}
The right adjoint of the stable splitting Theorem \ref{thm:bbcgsplitting} induces the product structure
in cohomology of a polyhedral product by observing the consequences of suspending the diagonal map.
We begin by defining 
$$\mathcal{H}^{q}(K;(\underline{X}, \underline{A})) = 
\bigoplus_{I\subseteq [m]}H^{q}\big(\widehat{Z}\big(K_I;(\underline{X}, \underline{A})_I\big).$$

\nd Next, we  impose an algebra structure on $\mathcal{H}^{q}(K;(\underline{X}, \underline{A}))$. In 
\cite[Section 3]{bbcg3}, it is shown that whenever  $I = J \cup L$,  there are well defined maps, 
{\em partial diagonals\/},
\newpage
$$\widehat{\Delta}^{J,L}_I \colon \widehat{Z}\big(K_I;(\underline{X}, \underline{A})_I\big) \longrightarrow
\widehat{Z}\big(K_J;(\underline{X}, \underline{A})_J\big) \wedge \widehat{Z}\big(K_L;(\underline{X}, \underline{A})_L\big)$$
making the diagram below commute
\begin{equation}\label{eqn:partial.diags}
\begin{CD}
Z(K;\xa)                       @>{\Delta_K}>>                     Z(K;\xa)\wedge Z(K;\xa)\\
@VV{\widehat{\Pi}_I}V                                               @V{\widehat{\Pi}_J\wedge\widehat{\Pi}_L}VV \\
\widehat{Z}\big(K_I;(\underline{X}, \underline{A})_I\big) @>{\widehat{\Delta}^{J,L}_I}>>
\widehat{Z}\big(K_J;(\underline{X}, \underline{A})_J\big)\wedge\widehat{Z}\big(K_L;(\underline{X}, \underline{A})_L\big),
\end{CD}
\end{equation}
\skp{0.2}
\nd where $\Delta_K$ is the reduced diagonal map. The maps $\Pi_I$, $\Pi_J$ and $\Pi_L$ are induced by
the appropriate projection maps
$$\Pi_S\colon Y^{[m]} = Y_1\times Y_2 \times\cdots\times Y_m \longrightarrow 
Y_{s_1} \wedge  Y_{s_2} \wedge \cdots \wedge Y_{s_k}$$
\nd where $S = \{s_1,s_2,\ldots,s_k\}$ corresponds in turn to $I$, $J$ and $L$. This allows us to define a product on
$\mathcal{H}^{q}(K;(\underline{X}, \underline{A}))$ as follows. Given cohomology classes 
\index{Polyhedral product, star product}
$u\in H^p\big(\widehat{Z}\big(K_J;(\underline{X}, \underline{A})_J\big)\big)$ and  
$v\in H^q\big(\widehat{Z}\big(K_L;(\underline{X}, \underline{A})_L\big)\big)$, define
 $$u\ast v=(\widehat{\Delta}^{J,L}_I)^\ast(u\otimes v)\; \in\; H^{p+q}(\widehat{Z}(K_I)).$$ 
The element $u\ast v\in H^{p+q}(\widehat{Z}(K_I))$ is called the
$\ast$-product of $u$ and $v$. Moreover, the commutativity of diagram \eqref{eqn:partial.diags} gives
\begin{equation}\label{eqn:star.cup} \
\widehat{\Pi}^{\ast}_I(u\ast v)= \widehat{\Pi}^{\ast}_J(u)\cup
\widehat{\Pi}^*_L(v) 
\end{equation} 
\nd where $\cup$ is the cup product for the CW-complex $Z(K;\xa)$.
 So, the $\ast$-product gives 
$\mathcal{H}^{q}(K;(\underline{X}, \underline{A}))$ a ring structure. 
Next, consider  the map
$$\eta=\bigoplus_{I\subseteq [m]}\Pi^{\ast}_I\colon\; \mathcal H^{\ast}(K;(\underline{X},\underline{A}))
\longrightarrow H^{\ast}\big({Z}(K;(\underline{X},\underline{A}))\big),$$
\nd which the stable splitting \eqref{eqn:splitting} ensures is an additive isomorphism. We can now state the main 
theorem of this section.
\begin{thm}[\cite{bbcg3}]\label{thm:products}
Let $K$ be an abstract simplicial complex with $m$ vertices and assume
that $(\underline{X},\underline{A})=\{(X_i,A_i, x_i)\}^m_{i=1}$ is a
family of based CW-pairs. Then
$$\eta:\mathcal H^*(K; (\underline{X},\underline{A}))\to H^*\big(Z(K; (\underline{X},\underline{A}))\big)$$ is a ring isomorphism.
\end{thm}
\nd  It follows that the $\ast$-product gives 
$\mathcal{H}^{q}(K;(\underline{X}, \underline{A}))$ an algebra structure. More details can be found in \cite[Section 3]{bbcg3}.
 This theorem generalizes the result of V.~Baskakov \cite{bask} \index{Baskakov's Theorem}
who proved it for the case of moment--angle complexes,
beginning with the cohomological splitting of Theorem \ref{thm:baskakov}. 

\begin{rem}
As discussed in \cite{bbcg3}, various corollaries follow from Theorem \ref{thm:products}.  Among these are to be found:
 \begin{enumerate}\itemsep2mm
\item If $\xa$ consists of suspension pairs and $J \cap L \neq \varnothing$, then $u\ast v = 0$ for all
$u\in H^p\big(\widehat{Z}\big(K_J;(\underline{X}, \underline{A})_J\big)\big)$ and  
$v\in H^q\big(\widehat{Z}\big(K_L;(\underline{X}, \underline{A})_L\big)\big)$.
\item Whenever $(\underline{CX},\underline{X})=\{(CX_i,X_i,x_i)\}^m_{i=1}$ is such that any 
finite product of the $X_i$ with the spaces $Z(K_I;(D^1,S^0))$ satisfies the strong form of the
K\"unneth theorem, the cup product structure for the
cohomology algebra $H^*(Z(K;(\underline{CX},\underline{X})))$ is a
functor of the cohomology algebras of $X$, and $Z(K_I;(D^1,S^0)))$
for all $I$. Real moment--angle complexes $Z(L;(D^1,S^0)))$ have been studied extensively in the work of 
A.~Suciu and A.~Trevisan \cite{suciu}, 
A.~Trevisan \cite{at}, L.~Cai \cite{cai1} and L.~Cai and S.~Choi \cite{caichoi}. 
\end{enumerate} 
\end{rem}

 We turn our attention now to the computation of the cohomology groups of a polyhedral product 
$H^\ast\big(Z(K;\xa)\big)$. \index{polyhedral product, spectral sequence}
For a suitable collection of pairs of spaces $(\underline{X}, \underline{A})$,  the cohomology of the polyhedral product 
$Z\big(K;  (\underline{X}, \underline{A}) \big)$  is computed using a spectral sequence 
 by the authors in \cite{bbcg10}.   A computation using different methods by Q.~Zheng can be found in 
\cite{zheng}.  The family of pairs
$(\underline{X}, \underline{A})$
 amenable to a cohomology calculation of  $Z\big(K;  (\underline{X}, \underline{A}) \big)$ satisfy a {\it strong freeness condition\/},
 \cite[Definition $2.2$]{bbcg10}.
 \begin{defin}\label{defn:strongfc}
The homology of $\xa$ with coefficients in a ring $k$, is said to be {\em strongly free\/} if the long exact sequence
\begin{equation}\label{eq:les}
\overset{\delta}{\to} \widetilde{H}^*(X_i/A_i) \overset{\ell}{\to} \widetilde{H}^*(X_i) \overset{\iota}{\to} \widetilde{H}^*(A_i)  
\overset{\delta}{\to} \widetilde{H}^{*+1}(X_i/A_i) \to  \end{equation}
\skp{0.1}
\nd satisfies the condition that there exist isomorphisms 
\nd \skp{0.2}
\begin{enumerate}\itemsep2mm
\item $\widetilde{H}^{\ast}(A_i)\cong E_i\oplus B_i$.
\item  $\widetilde{H}^{\ast}(X_i)\cong B_i \oplus C_i$,  \; where $B_i 
\underset{\simeq}{\overset{\iota}{\to}} B_i, \;\;\left.\iota\right|_{C_i} = 0$ 
\item  $\widetilde{H}^{\ast}(X_i/A_i)\cong C_i \oplus W_i$, \;
where $C_i \underset{\simeq}{\overset{\ell}{\to}} C_i, \;\; \left.\ell\right|_{W_i} = 0, \;\; 
E_i \underset{\simeq}{\overset{\delta}{\to}} W_i$ 
\end{enumerate}
for graded modules  $E_i, B_i, C_i$ and $W_i$ of finite type and free over $k$. 
\end{defin}
	
\nd  In particular,  for finite dimensional complexes,  this happens when coefficients are taken in a field.  
In order to describe the cohomology, more notation is required.   
For a simplicial complex 
$K$ with vertices in $[m] = \{1,2,3,\ldots, m\}$, and a subset $J=\{j_1,\cdots,j_r\}\subset [m]$,  
set $E^J = E_{j_1} \otimes \cdots \otimes E_{j_r}$.  Let $K_J$   be the full subcomplex with 
vertices in $J$. 
\nd For  $I=\{i_1,\cdots, i_t\} \subset [m]$,  and  for $\sigma \subset I$, set
$Y^{I,\sigma} = Y_1 \otimes \cdots  \otimes Y_t$ where 
	 $$ Y_j=  \begin{cases}	
	C_{i_j} & \mbox{ if }  i_j \in \sigma \\ 
	B_{i_j} & \mbox{ if } i_j \notin \sigma.
	\end{cases}$$	
	 and for $I=\varnothing$, set $Y^{\varnothing,\varnothing} = k$. 
\begin{thm}\label{thm:cohomology}
If $\xa$ satisfies the strong freeness condition there is a direct sum decomposition of the cohomology group
$$\widetilde{H}^{\ast}\big(\widehat{Z}\color{black}(K;\xa)\big) = \underset{I\cup 
J = [m], \;I \cap J =\varnothing}{\underset{\sigma \subset I, \;\sigma \in K}
		{\bigoplus}}E^J \otimes \widetilde{H}^{\ast}\big(\Sigma \big|{\rm{lk}}_{\sigma}(K_J)\big|\big) \otimes Y^{I,\sigma}$$
\nd with the convention that $\widetilde{H}^{\ast}(\Sigma \varnothing) = k$. \text{(Recall that ${\rm{lk}}_{\sigma}(K_J)$ is 
defined in Definition \ref{defn:link}.)}
\end{thm}
\nd Combined with the splitting Theorem \ref{thm:bbcgsplitting} in Section \ref{sec:splitting}, this theorem gives a complete
description of the cohomology of the topological spaces $Z\big(K;(\underline{X}, \underline{A})\big)$ with appropriate
coefficients. Moreover, the theorem generalizes 
to any multiplicative cohomology theory $h^*$ for which $\xa$ satisfies the 
natural formulation of the strong freeness condition for $h^{\ast}$. 
More information about the cohomology ring may be 
found in \cite{bbcg10} and \cite{zheng}.   The use of the methods developed in \cite{bbcg10} is illustrated next 
with a few examples, each having a different flavor. 

\begin{exm}\label{exm:sr}
In the particular case that  $H^\ast$ is cohomology with coefficients in a field $k$,  (or over $\mathbb{Z}$, under suitable 
freeness 
conditions), and the map $H^\ast(X_i)  \longrightarrow H^\ast(A_i)$ is surjective for all $i = 1,2,\ldots,m$, the spectral 
sequence constructed in \cite{bbcg10} collapses by
\cite[Proposition 3.8]{bbcg10}, and allows for a more concise statement of Theorem  \ref{thm:cohomology}. 
The assumption allows us to consider  $H^\ast(X_i/A_i)$ as a subring of 
$H^\ast(X_i)$ \index{Stannley--Reisner ring, generalized} and then
\cite[Corollary 3.9]{bbcg10} gives an isomorphism of rings
$$H^\ast\big(Z\big(K;(\underline{X}, \underline{A})\big)\big) \cong 
H^\ast(X_1) \otimes H^\ast(X_2) \otimes \cdots \otimes H^\ast(X_m)\big/I$$
\skp{0.1}
\nd where $I$ is the ideal generated by $H^\ast(X_{j_1}\big/A_{j_1}) \otimes H^\ast(X_{j_2}\big/A_{j_2}) 
\otimes \cdots \otimes H^\ast(X_{j_t}\big/A_{j_t})$ 
with $(j_1,j_2,\ldots,j_t)$ not spanning a simplex in $K$.
\end{exm}
\begin{rem}
Notice that this generalizes Corollary \ref{cor:sr}:
$$H^\ast\big(Z\big(K;(\mathbb{C}P^\infty, \ast);k\big)\big)\; \cong \;k(K)$$
\nd where $k(K)$ denotes the Stanley--Reisner ring of the simplicial complex $K$.
\end{rem}
\begin{exm}\label{exm:ca}
Another tractable example is the important case 
$$\xa =  (\underline{CA},\underline{A})$$
\nd  which includes the  example of  moment--angle complexes. In this case, with  $H^{\ast}(A_i)$  
free, the modules in the strong freeness condition of Theorem \ref{thm:cohomology}  
are 
$$E_i=\widetilde{H}^\ast(A_i), \quad B_i=0,\quad  C_i=0, \;\text{and}\; W_i =\widetilde{H}^{\ast}(\Sigma A_i).$$ 
The modules  $Y^{I,\sigma}$ of Theorem \ref{thm:cohomology}, take a particularly simple form:
$ Y^{I,\sigma}= 0$ if $I \neq \varnothing$  and is equal to $k$  if  $I=\varnothing$, (see \eqref{eq:les}, item (1)).
As a consequence, the only links that appear in Theorem \ref{thm:cohomology}  are $\text{lk}_{\varnothing}(K)= K$ and 
we have 
\begin{equation}\label{eq:cx}
\widetilde{H}^\ast\big(\widehat{Z}(K;(\underline{CA},\underline{A})\big) =  \widetilde{H}^\ast(\Sigma \big|K\big|) 
\otimes \widetilde{H}^\ast(A_{1}) \otimes \cdots \otimes \widetilde{H}^\ast(A_{m}).  \end{equation}
\nd This result agrees with the one given by the wedge lemma \cite{wzz} as described in \cite{bbcg2}.  Finally, the splitting 
theorem
\cite[Theorem $2.10$]{bbcg2} gives the cohomology of the polyhedral product as
\begin{align*}
 H^{\ast}\big(Z(K;(\underline{CA},\underline{A}))\big)\;&\;= \underset{I \subset [m]}{\bigoplus} H^{\ast}(\Sigma|K_I|) 
\otimes \widetilde{H}^{\ast}(A_{i_1}) \otimes \cdots \otimes \widetilde{H}^{\ast}(A_{i_t})\\
&\;\subset \; H^{\ast}\big(Z(K;(D^1,S^0))\big) \otimes H^{\ast}(A_{1}) \otimes \cdots \otimes H^{\ast}(A_{m}).
 \end{align*}\itemsep3mm
 \skp{0.1}
 \nd There is now an evident product on $H^\ast\big(Z(K; (\underline{CA},\underline{A}))\big)$ induced by coordinate-wise 
 multiplication
 and by the product in $H^\ast\big(Z(K; (D^1,S^0))\big)$, 
which is the case corresponding to $A_i = S^0$. The latter is described  in \cite{cai1}. 
In \cite{bbcg10}  it is shown, using the results of \cite{bbcg3}, that  this is indeed the ring structure in 
$H^\ast\big(Z(K;(\underline{CA},\underline{A}))\big)$. 
\end{exm}
\begin{rem} The ring structure for general $\xa$ satisfying the strong freeness condition 
may be found in \cite{bbcg3}, \cite{bbcg10} and \cite{zheng}.\end{rem}

The next example is of a different nature.
 \begin{exm}\label{exm:cohomology}  
 Consider a CW-pair $(X,A)$ with cohomology satisfying
 \begin{equation}\label{eqn:example}
 H^{\ast}(X) =\mathbb{Z}\{ b_4, c_6\}  \quad \text{and} \quad H^{\ast}(A)=\mathbb{Z}\{e_2, b_4\}
 \end{equation}
 \nd where the dimensions of the classes are given by the subscripts. A trivial example of $\xa$ is given by wedges of spheres in 
 the appropriate dimensions. A more interesting example is obtained from the mapping cylinder of the composite map
 \begin{equation}\label{eqn:cp2cp3}
 f\colon \mathbb{C}P^2 \hookrightarrow \mathbb{C}P^3 \to \mathbb{C}P^3/\mathbb{C}P^1.
 \end{equation}
 \nd We denote the mapping cylinder of \eqref{eqn:cp2cp3} by $M_f$ and consider the pair
 $(M_f, \mathbb{C}P^2)$, which satisfies
 the cohomology condition above.
 \begin{rem} Notice here that Theorem \ref{thm:cohomology} will give the {\em same} cohomology for $\widehat{Z}(K ; \xa)$ 
 whether we realize condition \eqref{eqn:example} with a pair made from appropriate wedges of spheres, or from the 
 projective spaces above.  Given $K$, the cohomology of $\widehat{Z}(K ; \xa)$ depends on the modules 
 $E_i, B_i$ and $C_i$ only. 
 \end{rem} 
 
 For the example at hand, let $K$ be the simplicial complex with 3 vertices and  edges $\{1,3\}, \{1,2\}$.
 The cases in this example are indexed by the $I$ in Theorem \ref{thm:cohomology} starting with $I = \varnothing$ and building up 
to $I=\{1,2,3\}$.  For each $I$ there are the sub-cases indexed by the simplices $\sigma \subset I$.  Theorem \ref{thm:cohomology}
reduces the calculation to bookkeeping, as it does with every example, but the bookkeeping can become quite complicated. 
\begin{enumerate}\itemsep1.5mm 	 
\item $I=\varnothing$, $\sigma = \varnothing$.\\
Here $J = \{1,2,3\}$ and  $|$lk$_{\varnothing}(K_J)| = |K|$ is contractible. So, there
is no contribution to the Poincar\'{e} series for $H^\ast\big(\widehat{Z}(K; \xa)\big)$ in this case.
\item $I=\{1\}$, $\sigma = \varnothing$.\\			
Now, $J  =\{2,3\}$, so $E^J$ contributes $(t^2)^2$, $Y^{I,\varnothing}= b_4$ which gives a $t^4$, 
and $|$lk$_{\varnothing}(K_J)| = |\{\{2\},\{3\}\}|=S^0$. So for this case, Theorem \ref{thm:cohomology} specifies a 
total contribution of $t^9$ to the Poincar\'{e} series for $H^\ast\big(\widehat{Z}(K; \xa)\big)$.
\item  $I=\{1\}$, $\sigma=\{1\}$.\\
Again, $J  =\{2,3\}$, so $E^J$ contributes $(t^2)^2$, $Y^{I,\{1\}}= c_6$ which gives a $t^6$, 
and $|$lk$_{\{1\}}(K_J)| = |\{\{2\},\{3\}\}|=S^0$. So for this case, Theorem \ref{thm:cohomology} specifies a 
total contribution of $t^{11}$ to the Poincar\'{e} series for $H^\ast\big(\widehat{Z}(K; \xa)\big)$.
\end{enumerate} 
\nd Continuing in this way, we arrive at the (reduced) Poincar\'{e} series
$$\overline{P}\big(H^*(\widehat{Z}(K; \xa))\big) = t^9+t^{11}+3t^{12}+5t^{14}+2t^{16}.$$
\end{exm}	
\begin{rem}
 This illustrative example  lends itself to direct calculation.
For $K$ as above, it follows directly from the definition that 
$$\widehat{Z}(K;(X,A)) = X\wedge \big((X\wedge A)\cup(A\wedge X)\big).$$ 
\nd The K\"unneth theorem now reduces the calculation
to $H^\ast\big((X\wedge A)\cup(A\wedge X)\big)$, which is direct for the pair $(X,A)$ in Example \ref{exm:cohomology} 
by the Mayer-Vietoris sequence.
\end{rem}
\skp{0.2}		

\section{A geometric approach to the cohomology of polyhedral products}\label{sec:cartan}
In the forthcoming paper \cite{bbcg11}, the authors show that for certain pairs $\xa$, called {\em wedge decomposable\/},
the algebraic decomposition given by Theorem \ref{thm:cohomology} is a consequence of an underlying geometric splitting. 
{\em Moreover, the consequences of this observation extend to more general based CW-pairs\/}. (The results of 
this section are from the authors' unpublished preprint from 2014, which in turn originated from an earlier preprint from 2010,
and is currently being revised.)

\begin{defin}\label{def:wedgedecomp} \index{Wedge decomposable pair}
Based CW-pairs  of the form 
$$(\underline{X}, \underline{A}) = (\underline{B\vee C}, \underline{B\vee E})$$
so that, for all $i$, $(X_{i},A_{i }) = (B_{i}\vee C_{i}, B_{i}\vee E_{i})$, where $E_{i}\hookrightarrow C_{i} $ is a
null homotopic inclusion, are called wedge decomposable.
\end{defin} 
For such pairs, there is a decomposition of  the smash polyhedral product.  Let 
$$J = \{j_1,j_2,\ldots,j_k\} \subset [m]$$
and set $\widehat{B}^J = B_{j_1}\wedge B_{j_2}\wedge \cdots \wedge B_{j_k}$. Similarly, define $\widehat{C}^J$. 
\begin{thm}[\cite{bbcg11}]\label{thm:cartan}
Let $(\underline{X}, \underline{A}) = (\underline{B\vee C}, \underline{B\vee E})$ be a wedge decomposable
pair. Then there is a homotopy equivalence
$$\widehat{Z}\big(K;(\underline{X}, \underline{A})\big) \longrightarrow 
\bigvee_{I\leq [m]}\widehat{Z}\big(K_{I};(\underline{C},\underline{E})_I\big)\wedge \widehat{B}^{([m]-I)}$$
\nd which is natural with respect to maps of wedge decomposable pairs.
\end{thm} 
\nd Since  the inclusion $E_{i}\hookrightarrow C_{i}$ is null homotopic, the terms 
$\widehat{Z}\big(K_{I};(\underline{C},\underline{E})_I\big)$ are completely determined by the 
the {\em wedge lemma\/}, \cite{bbcg1}, as follows.
\begin{prop}\label{prop:wedge}
For pairs $(\underline{C},\underline{E})$ as above there is a homotopy equivalence
 $$\widehat{Z}\big(K_{I};(\underline{C},\underline{E})_I\big) \to   
\bigvee_{\sigma \in K_{I}} |{\rm{lk}}_{\sigma}(K_{I})|\ast \widehat{D}_{\underline{C},\underline{E}}^{I}(\sigma) $$    
\nd where $ |{\rm{lk}}_{\sigma}(K_{I})|$ is  the realization of the link of $\sigma$ in the
full subcomplex $K_{I}$ and 
\begin{equation}\label{eqn:d.sigmace}
\widehat{D}_{\underline{C},\underline{E}}^{I}(\sigma)  =\bigwedge^{|I|}_{j=1}W_{i_{j}},\quad {\rm with}\quad
W_{i_{j}}=\left\{\begin{array}{lcl}
C_{i_{j}}  &{\rm if} & i_{j}\in \sigma\\
E_{i_{j}}  &{\rm if} & i_{j}\in I-\sigma.
\end{array}\right.
\end{equation}
\end{prop}
\nd Combined with the splitting Theorem \ref{thm:bbcgsplitting}, these results give a complete
description of the topological spaces $Z\big(K;(\underline{X}, \underline{A})\big)$ for wedge decomposable pairs $\xa$.

The case $E_i \simeq \ast$ simplifies further by Theorem \ref{thm:more.contractible.A} to give the next corollary.
\begin{cor}\label{cor:E.a.point}
For wedge decomposable pairs of the form $(\underline{B\vee C}, \underline{B})$, corresponding to $E_i \simeq \ast$
for all $i = 1,2,\ldots,m$, there are homotopy equivalences
$$\widehat{Z}\big(K_{I};(\underline{C},\underline{E})_I\big) \simeq \widehat{Z}\big(K_{I};(\underline{C},\ast)_I\big)
\simeq  \widehat{C}^I,$$
and so Theorem \ref{thm:cartan} gives $\widehat{Z}\big(K;(\underline{B\vee C}, \underline{B})\big) \simeq
\bigvee_{I\subseteq [m]}\big(\widehat{C}^I\wedge \widehat{B}^{([m]-I)}\big)$.
\end{cor}

\nd Notice that the Poincar\'e series for the space $\widehat{Z}\big(K;(\underline{B\vee C}, \underline{B})\big)$
follows easily from Corollary \ref{cor:E.a.point}. 
\begin{rem} In comparing these observations with  Theorem \ref{thm:cohomology}, notice that the links appear in the terms
$\widehat{Z}\big(K_{I};(\underline{C},\underline{E})_I\big)$. Also, while  Theorem \ref{thm:cartan} and Proposition \ref{prop:wedge} 
give a geometric underpinning for the cohomology calculation in Theorem \ref{thm:cohomology} for wedge decomposable pairs
only, the geometric splitting does not require that $E,B$ or $C$ have torsion-free cohomology
\end{rem}
Theorem \ref{thm:cartan} applies particularly well in cases where spaces have unstable attaching maps.
\begin{exm}
The homotopy equivalence $S^{1}\wedge Y \simeq \Sigma(Y)$ implies homotopy equivalences
\begin{equation}\label{eqn:he1}
\Sigma^{mq}\big(\widehat{Z}(K;(\underline{X},\underline{A}))\big)
\longrightarrow \widehat{Z}\big(K;\big(\underline{\Sigma^{q}(X)},\underline{\Sigma^{q}(A)}\big)\big)
\end{equation}
\nd where as usual,  $m$  is the number of vertices of $K$. Next, recalling that $SO(3) \cong \mathbb{R}\rm{P}^{3}$,
consider the pair
$$(X,A) = \big(SO(3), \mathbb{R}\rm{P}^{2}\big)$$
\nd for which there is a well known homotopy equivalence of pairs
\begin{equation}\label{eqn:he2}
\big(\Sigma^{2}\big(SO(3)\big), \;\Sigma^{2}\big(\mathbb{R}\rm{P}^{2}\big)\big) \longrightarrow 
\big(\Sigma^{2}\big(\mathbb{R}\rm{P}^{2}\big)\vee \Sigma^{2}(S^{3}),\;
\Sigma^{2}\big(\mathbb{R}\rm{P}^{2}\big)\big),
\end{equation}
\nd making the pair $\big(SO(3), \mathbb{R}\rm{P}^{2}\big)$ {\em stably wedge decomposable\/}.
\nd Now, combining \eqref{eqn:he1} and \eqref{eqn:he2}, we get a homotopy equivalence
$$\Sigma^{2m}\big(\widehat{Z}(K;(SO(3),\mathbb{R}\rm{P}^{2}))\big)
\longrightarrow 
\widehat{Z}\big(K;\big(\Sigma^{2}(\mathbb{R}\rm{P}^{2})\vee \Sigma^{2}(S^{3}),\;
\Sigma^{2}(\mathbb{R}\rm{P}^{2})\big).$$
\nd Finally, Theorem \ref{thm:cartan} allows us to conclude that $\widehat{Z}(K;(SO(3),\mathbb{R}\rm{P}^{2}))\big)$,
and hence the polyhedral product $Z(K;(SO(3),\mathbb{R}\rm{P}^{2}))$,
is stably a wedge of smash products of $S^{3}$ and $\mathbb{R}\rm{P}^{2}$. Similar splitting results follow for
the polyhedral product whenever the spaces $X$ and $A$ split after finitely many suspensions. In particular, the fact that
$\Omega^{2}S^{3}$ splits stably into a wedge of Brown--Gitler spectra implies that the polyhedral product 
$Z\big(K; (\Omega^{2}S^{3}, \ast)\big)$ splits  stably  into a wedge of smash products of Brown--Gitler spectra.
\end{exm}

The result of the previous section can be exploited to give information about 
the groups $H^{*}\big(\widehat{Z}(K;(\underline{X},\underline{A}))\big)$ over a field $k$ for pointed, 
finite, path connected  finite pairs of 
CW-complexes $(\underline{X}, \underline{A})$, which are {\bf {\em not\/}} wedge decomposable. In so doing, we 
explain further the remark in Example \ref{exm:cohomology}.

Given $(\underline{X}, \underline{A})$, let $B_i$, $C_i$ and $E_i$ be the $k$--modules  specified in items (1) and (2) of 
Definition \ref{defn:strongfc}. 
Now, wedges of spheres and Moore spaces $B'_i$ , $C'_i$ and $E'_i$ exist realizing the modules $B_i$, $C_i$ and $E_i$ so that
$$ (\underline{B'\vee C'},\underline{B'\vee E'})$$ 
\nd satisfies the criterion for a wedge decomposable pair as in Definition \ref{def:wedgedecomp}. Moreover, the diagram
below commutes:
{\begin{equation*}
\begin{CD}
H^{\ast}(B'_{j}\vee C'_{j})  @>>{\cong}> H^{\ast}(X_{j})\\
@VV{}V           @VV{}V          \\
H^{\ast}(B'_{j}\vee E'_{j})@>{\cong}>>  H^{\ast}(A_{j}).
\end{CD}
\nd \end{equation*}}
\nd This leads to the following result.

\begin{thm}[\cite{bbcg11}]\label{thm:cartan2}
Under the conditions stated above, the following isomorphism of groups holds for cohomology with coefficients
in a field $k$ 
$$H^{*}\big(\widehat{Z}(K;(\underline{X},\underline{A}))\big) \cong
 H^{*}\big(\widehat{Z}(K;(\underline{B'\vee C'},\underline{B'\vee E'}))\big)$$ 
\nd where the right hand side is determined by Theorem \ref{thm:cartan} and Proposition \ref{prop:wedge}.
 (This is not necessarily an isomorphism of modules over the Steenrod algebra.) 
\end{thm}
As a consequence, the additive structure of $H^{\ast} \big(Z(K;(X,A));k\big)$ over a 
field $k$, is a functor which is determined by   $H^{\ast}(X)$, $H^{\ast}(A)$, and the ranks of the 
restriction maps appearing in the extension 
$$0 \to V_i \to H^i(X) \to  H^i(A) \to W_i \to 0, $$ 
(cf. \cite{zheng}).  Information about the ring structure requires additional assumptions. 

\section{Polyhedral products and the Golodness of  monomial ideal rings}\label{sec:golodness}
We begin with the definition of a Golod ring, \cite{frank, lk}. \index{Golodness of rings}
\begin{defin}
Let $S = k[x_1,x_2,\ldots,x_m]$ be a polynomial ring in $m$ variables over a field $k$, and let $I = (m_1,m_2,\ldots,m_r)$
be an ideal generated by monomials. The monomial ideal ring $R = S/I$ is called {\em Golod\/} if
\begin{equation}\label{eqn:golod}
\sum_{j=0}^{\infty}{\rm dim}\hspace{0.5mm}{\rm Tor}_j^{R}(k,k)t^j\;\; = 
\;\;\frac{(1+t)^m}{1 - t(\ds{\sum_{j=0}^{\infty}{\rm dim}\hspace{0.5mm}{\rm Tor}_j^{S}(R,k)t^j-1)}}.
\end{equation}
\end{defin}
 \begin{rem}
 According to \cite{lk} and \cite{golod}, Serre had observed that  the coefficients on the left hand side 
in \eqref{eqn:golod} are always less than or equal to the corresponding coefficients on the right. 
\end{rem}
 In \cite{golod}, Golod showed that \eqref{eqn:golod}
holds if and only if all products and all higher Massey products vanish in To$r^{S}(R,K)$, (see also \cite{lk}).  
Higher Massey products
in $H^\ast(X;\mathbb{Q})$ obstruct the rational formality of $X$. As is noted in \cite{limon3}, this has a bearing 
on whether 
a complex manifold can admit a K\"{a}hler structure, an observation about the formality of K\"{a}hler 
manifolds  which can be found  in the paper by  P.~Deligne, Ph.~ A.~Griffiths, J.~W.~Morgan, and 
D.~Sullivan \cite{dgms}.  Golodness is relevant also in symplectic geometry and the theory of subspace 
arrangements; an 
overview of these connections can be found in \cite{ds} and\cite{limon3}.

The study of Golodness is particularly important
in the case that $R$ is the Stanley--Reisner ring $k(K)$ of a simplicial complex $K$,  (see \eqref{eqn:srring}). 
The property has been 
much investigated by algebraic combinatorial theorists. The arrival of moment--angle complexes has invigorated this line of 
research.  Moment--angle complexes  become relevant via the split  fibration  arising from 
Theorem \ref{eqn:djequalsbp}, namely
\begin{equation}\label{eqn:loopsmac}
\Omega{Z}(K;(D^2,S^1)) \longrightarrow \Omega{Z}\big(K;(\mathbb{C}{\rm P}^\infty, \ast)) \longrightarrow T^m.
\end{equation}
As noted in Section \ref{sec:cohomofamac},  there is an isomorphism  of rings
\begin{equation}\label{eqn:maccohom} 
H^{\ast}\big(Z(K;(D^2,S^1))\big) \cong {\rm{Tor}}^*_{S}(k(K),k).
\end{equation}
Using these ideas, and the isomorphism
$$H^{\ast}\big( \Omega{Z}\big(K;(\mathbb{C}{\rm P}^\infty, \ast)\big) \cong {\rm{Tor}}^{\ast}_{k(K)}(k,k),$$
due to \nd Buchstaber and Panov \cite{bp1}, in conjunction with an Eilenberg--Moore spectral sequence argument, 
J.~Grbi\'c  and S.~Theriault, \cite[Theorem $11.1$]{gtshifted}, 
were able to recapture for Stanley--Reisner rings the inequality involving the terms in \eqref{eqn:golod}, mentioned in the
remark above
and attributed to Serre in the general case. Moreover, extending Theorem \ref{thm:gtshifted}, they showed that
for a certain class $\mathcal{F}_0$, of simplicial complexes, obtained from shifted simplicial complexes by elementary 
topological operations, $Z(K;(D^2,S^1))$ has the homotopy type of a wedge of 
spheres, implying Golodness. G.~Denham and A.~Suciu also study obstructions to Golodness, and in particular, triple Massey 
products in lowest degree in \cite{ds}.
\begin{thm}\cite[Theorem $1.2$]{gtshifted}.
If $K\in \mathcal{F}_0$, then $k(K)$ is a Golod ring over $k$.
\end{thm}
K.~Iriye and D.~Kishimoto \cite{ik5}, begin with the stable decomposition of Theorem \ref{thm:xicontractible} applied to the case
$(\underline{CX},\underline{X})$ to get 
\begin{equation}\label{eqn:stablecxx}
\Sigma{Z}(K;(\underline{CX},\underline{X}))\; \simeq\; 
\Sigma\!\!\!\bigvee_{\varnothing \neq I \subset [m]}\!\!\Sigma{K_{I}}\wedge \widehat{X}^I.
\end{equation}
They notice that if \eqref{eqn:stablecxx} desuspends for $X = S^1$,  the space $Z(K;(D^2,S^1))$
remains a suspension and hence \eqref{eqn:maccohom} implies that $k(K)$ is Golod over $k$, \cite[Corollary $3.11$]{porter2}.

They introduce a structure on $Z(K;(D^1,S^0))$ which they call the {\em fat wedge filtration\/}.  \index{Fat wedge filtration}
This then allows them to give sufficient conditions on $K$ ensuring  the desuspension of 
\eqref{eqn:stablecxx} for any $X$ and  concluding that $k(K)$ is a Golod ring. Among these results is the next theorem.
\index{Polyhedral product, wedge of spheres}

\begin{thm}\cite{ik5}\label{thm:shellable} 
If the Alexander dual of $K $is shellable, then $Z(K;(D^n,S^{n-1}))$ has the homotopy
type of a wedge of spheres, and hence $k(K)$ is Golod over $k$.
\end{thm}

They are able also to recover a result of Herzog, Reiner and Welker.
\begin{thm}\cite{hrw, ik5}
If the Alexander dual of the simplicial complex $K$ is sequentially Cohen--Macaulay over $k$, then $k(K)$ is Golod over
$k$. 
\end{thm}
Fat wedge filtration techniques are also employed by K.~Iriye and D.~Kishimoto to prove that if $K$ is a triangulated surface
orientable over $k$, then $k(K)$ is Golod over $k$ if and only if it is {\em 2--neighborly\/}, that is, any two vertices are 
connected by an edge, \cite[Theorem $1.3$]{ik3}. Using similar ideas,  K.~Iriye and D.~Kishimoto \cite{ik4} find a two-dimensional 
simplicial complex $K$ with $k(K)$ Golod over $\mathbb{Q}$ but not over $\mathbb{Z}/p$.
\begin{rem}
The fat wedge filtration methods of K.~Iriye and D.~Kishimoto \cite{ik5} allow them to recover a variety of  stable and unstable
homotopy decompositions of moment--angle complexes and other polyhedral products. 
\end{rem}
K.~Iriye and T.~Yano \cite{iy}, also take the approach of desuspending  \eqref{eqn:stablecxx}. They use a simplicial complex 
constructed from a triangulated Hopf map $\eta\colon S^3 \to S^2$, and Alexander duality to prove the existence of
a simplicial complex $K$ for which $k(K)$ is a Golod ring but $Z(K;(D^2,S^1))$ is not a suspension. 

J.~Grbi\'c, T.~Panov, S.~Theriault and J.~Wu give an example, \cite[Example $3.3$]{gptw}, of a simplicial complex $K$, 
the standard 6-vertex triangulation of $\mathbb{R}{\rm P}^2$, which is Golod
over any field but $Z(K;(D^2,S^1))$ is not a wedge of spheres  because it has torsion. I.~Limonchenko gives a 
similar example arising from
a 9-vertex triangulation of $\mathbb{C}{\rm P}^2$, \cite[Theorem $2.5$]{limon2}. He does this by explicitly computing
the Betti numbers of $Z(K;(D^2,S^1))$ using \eqref{eqn:betti}, and then checking that the attaching maps in the
stable splitting Corollary \ref{cor:bigmac} are all stable maps. In \cite{limon3}, I.~Limonchenko determines conditions 
on the bigraded Betti numbers \eqref{eqn:betti} which imply the non-Golodness of $K$ over $k$. He uses the simplicial
wedge construction \eqref{eqn:wedge}, to find a family of generalized moment--angle complexes, such that for any
$l,r \geq 2$, the family contains an $l$-connected manifold $M$ with a non-trivial $r$-fold Massey product in 
$H^{\ast}(M, \mathbb{Q})$.

Another approach to the Golodness problem for $K$ via the homotopy theory of the moment--angle complex
$Z(K;(D^2,S^1))$, can be found in the work of 
P.~Beben and J.~Grbi\'c, \cite{bg1}. A.~Berglund gives a combinatorial condition in \cite{berglund} which ensures the 
Golodness of a monomial ideal ring.
\vspace{-1.4\baselineskip}
\section{Higher Whitehead products and loop spaces}\label{sec:hwp}
Higher Whitehead products were introduced into the homotopy theory of moment--angle complexes in the work
of T.~Panov and N.~Ray \cite{pr}. They were concerned with the problem of realizing sphere wedge summands
in $Z(K;(D^2,S^1))$ by {\em higher\/} Whitehead products. This work is summarized in 
\cite[Sections $8.4$ and $8.5$]{bp2}. Cases of this problem are discussed  also in  \cite{gt3} and \cite{ik6}. Recently,
S.~Abramyan \cite{abramyan} showed that not all sphere summands can be realized in this way using higher Whitehead
products.  These ideas are developed further in \cite{ap}. 

In the language of moment--angle complexes, the maps of \eqref{eqn:wp} can be reformulated as follows. Let
$K$ be the simplicial complex consisting of two  discrete  vertices,  the boundary of a one-simplex 
$\partial\Delta^1$,  and consider 
the corresponding  
fibration \eqref{eqn:djequalsbp},
\begin{equation}\label{eqn:basicfib2points}
Z(\partial\Delta^1;(D^2,S^1)) \stackrel{\widetilde{\omega}}\longrightarrow Z(K;(\mathbb{C}{\rm P}^\infty, \ast)) \longrightarrow 
\mathbb{C}{\rm P}^\infty \times \mathbb{C}{\rm P}^\infty.
\end{equation}
 We have 
\begin{align*}
Z(\partial{\Delta^1};(D^2,S^1)) &= S^1 \times D^2 \cup_{S^1\times S^1}D^2 \times S^1\; \simeq\; S^3\\
Z(\partial{\Delta^1};(S^2,\ast)) &= (S^2\times \ast)\cup_{\ast\times \ast}(\ast \times S^2) \;=\; S^2\vee S^2\\
Z(\partial{\Delta^1};(\mathbb{C}{\rm P}^\infty, \ast)) &= \mathbb{C}{\rm P}^\infty \vee \mathbb{C}{\rm P}^\infty\\
Z(\Delta^1;(\mathbb{C}{\rm P}^\infty, \ast)) &= \mathbb{C}{\rm P}^\infty \times \mathbb{C}{\rm P}^\infty.
\end{align*}
Then \eqref{eqn:basicfib2points} factors as
$$S^3 \;\simeq\; Z(\partial{\Delta^1};(D^2,S^1)) \longrightarrow Z(\partial{\Delta^1};(S^2,\ast)) \longrightarrow 
Z(\partial{\Delta^1};(\mathbb{C}{\rm P}^\infty, \ast))
\longrightarrow \mathbb{C}{\rm P}^\infty \times \mathbb{C}{\rm P}^\infty.$$
\nd Here, the first map is induced by the map of pairs $(D^2,S^1) \to (S^2,\ast)$, the 
second by the inclusion $(S^2,\ast) \to (\mathbb{C}{\rm P}^\infty, \ast)$, and the third by the inclusion of
simplicial complexes $\partial{\Delta^1} \hookrightarrow \Delta^1$.  

Now suppose that  $\partial{\Delta^{k-1}}$ 
is a minimal missing face of a simplicial complex $K$, then we get the generalization 
\begin{multline}\label{eqn:higherwp}
S^{2k-1} \;\simeq\; Z(\partial{\Delta^{k-1}};(D^2,S^1)) \longrightarrow Z(\partial{\Delta^{k-1}};(S^2,\ast))\\ \longrightarrow 
Z(\partial{\Delta^{k-1}};(\mathbb{C}{\rm P}^\infty, \ast)) \longrightarrow Z(K;(\mathbb{C}{\rm P}^\infty, \ast)). 
\end{multline}
Here, the last map is induced by the inclusion $\partial{\Delta^{k-1}} \hookrightarrow K$ and 
$Z(\partial{\Delta^{k-1}};(\mathbb{C}{\rm P}^\infty, \ast))$ retracts off $Z(K;(\mathbb{C}{\rm P}^\infty, \ast))$ by
 \cite[Proposition $3.3.1$]{ds}. Following \cite[page 339]{bp2}, the group 
 $\pi_{2}\big(Z(K;(\mathbb{C}{\rm P}^\infty,\ast))\big) \cong \mathbb{Z}^m$
 has $m$ generators
 $$\widehat{\mu}_i \colon S^2 \to \mathbb{C}{\rm P}^\infty \xrightarrow{i} 
 \mathbb{C}{\rm P}^\infty\vee \mathbb{C}{\rm P}^\infty \vee
 \cdots \vee \mathbb{C}{\rm P}^\infty \to Z(K;(\mathbb{C}{\rm P}^\infty,\ast)),$$
 \nd where the second map is the inclusion of the $i^{{\rm th}}$ wedge summand and the last map is induced by the 
 inclusion of the zero-skeleton into $K$.  Next, labelling the vertices of the missing face $\partial{\Delta^{k-1}}$ by 
 $\{i_1,i_2,\ldots,i_k\}$, we call the  composite map \eqref{eqn:higherwp} the {\em $k$-fold higher 
 Whitehead product\/} \index{Whitehead product, higher} and  denote it by the symbol 
 $[\widehat{\mu}_{i_{1}},\widehat{\mu}_{i_{2}},\ldots,\widehat{\mu}_{i_{k}}]_w$.
 
These maps, and associated Samelson products, play an important role in the study of the homotopy theory
of  $\Omega{Z}(K;(D^2,S^1))$ and  $\Omega{Z}(K;(\mathbb{C}{\rm P}^\infty,\ast))$, particularly cases for which
$Z(K;(D^2,S^1))$ is homotopy equivalent to a wedge of spheres, see \cite[Section $8.4$]{bp2}. In the case of a 
flag complex $K$, T.~Panov and N.~Ray \cite{pr} compute the rational Pontrjagin ring of 
$\Omega{Z}(K;(\mathbb{C}{\rm P}^\infty,\ast))$ by introducing various algebraic and geometric models. Motivated in part
by this, and the work of S.~Papadima and A.~Suciu \cite{papsu}, N.~Dobrinskaya \cite{nd} also  addresses 
$\Omega{Z}(K;(X,\ast))$ from a different point of view, and relates the 
computation to diagonal arrangements.  (Recall from Section \ref{sec:splitting} that $Z(K;(D^2,S^1))$ is a deformation 
retract of $Z(K;(\mathbb{C}, \mathbb{C}^{\ast}))$, the complement of a complex coordinate arrangement 
determined by $K$.) 

K.~Iriye and D.~Kishimoto, \cite{ik6},  call a simplicial complex $K$ {\em totally fillable} 
\index{Simplicial complex, totally fillable} 
if each of its full subcomplexes
has the property that it becomes contractible if some of its minimal missing faces are added. They employ their theory 
of {\em fat wedge filtrations\/} \cite{ik5}, to show that if $K$ is totally fillable, the moment--angle complex
$Z(K;(D^2,S^1))$ decomposes as a wedge of spheres. For each such sphere $S^\alpha$, they recognize the map
$$S^\alpha \hookrightarrow Z(K;(D^2,S^1)) \stackrel{\widetilde{\omega}}\longrightarrow Z(K;(\mathbb{C}{\rm P}^\infty, \ast))$$
in terms of higher and iterated Whitehead product, \cite[Theorem $1.3$]{ik6}.

K.~Iriye and D.~Kishimoto then extend the discussion, as do J.~Grbi\'c and S.~Theriault in \cite{gt3}, to the more general 
 fibration 
\begin{equation}\label{eqn:coneloopfib} \index{Loop spaces}
Z\big(K;(C(\Omega{\underline{X})\big), \Omega{\underline{X}})) 
\stackrel{\widetilde{\omega}}\longrightarrow  Z(K;(\underline{X},\ast))} \longrightarrow  \prod_{i=1}^{m}X_i,
\end{equation}
(which was studied also by  G.~Porter, \cite{porter2}). In both papers,  the map $\widetilde{\omega}$ is described fully 
in particular cases. More comprehensive information about Whitehead products in the context of 
\eqref{eqn:coneloopfib}  is to be found in the paper by S.~Theriault \cite{th2}. The rational type of the fibre in 
\eqref{eqn:coneloopfib}
is studied also in the paper of Y.~F\'elix and D.~Tanr\'e \cite{ft}.

\subsection*{Acknowledgements} The authors would like to thank Haynes Miller, Peter Landweber, Taras Panov and 
Jelena Grbi\'c for their 
careful reading of the manuscript and for their many valuable suggestions. Our thanks also to Santiago L\'opez de Medrano, 
Daisuke Kishimoto,  Matthias Franz,  Alex Suciu, Alvise Trevisan, Mentor Stafa, Stephen Theriault, Kouyemon Iriye and 
Graham Denham, for  comments and recommendations which have led to a discernible improvement in the presentation.

\bibliographystyle{amsalpha}

\end{document}